\documentclass[a4paper,12pt]{amsart}  
\usepackage{amssymb}   
\usepackage{amsmath,amsthm,geometry,graphicx} 
\usepackage{amscd}
\usepackage{psfrag,xspace} 
\usepackage[all,dvips,cmtip]{xy} 

\usepackage{amsmath}
\usepackage{amssymb}

\usepackage{psfrag,graphicx} 
\usepackage{fullpage,url}
\usepackage{amsfonts} 
\usepackage{latexsym} 
\usepackage{epsfig} 
\usepackage{amsthm}
\usepackage[all]{xy}
\usepackage{mathrsfs}
\usepackage{amsthm,amsmath,amssymb,amsfonts, amscd,graphics, latexsym,
enumerate, stmaryrd,xspace,verbatim,epic,eepic}

\numberwithin{equation}{section} 

\newtheorem{thm}[equation]{Theorem}
\newtheorem{cor}[equation]{Corollary}
\newtheorem{lem}[equation]{Lemma}
\newtheorem{prop}[equation]{Proposition}
\theoremstyle{definition} 
\newtheorem{defn}[equation]{Definition}
\newtheorem{defi}[equation]{Definition}

\newtheorem{exmp}[equation]{Example}%
\newtheorem*{notn}{Notation}%
\newtheorem{construction}[equation]{Construction}%

\newtheorem{rem}[equation]{Remark}

{\nolinebreak $\Box$ \end{trivlist}}

\newcommand{\OO}{\mathcal O}
\newcommand{\GG}{\mathcal G}
\newcommand{\XX}{\mathcal X}
\newcommand{\UU}{\mathcal U}
\newcommand{\YY}{\mathcal Y}

\newcommand{\DD}{\mathcal D}

\newcommand{\TT}{\mathcal T}

\newcommand{\Gm}{{\mathbb G}_m}

\newcommand{\eps}{\varepsilon}
\newcommand{\fs}{\!\!\!\fatslash}
\DeclareMathOperator{\sm}{sm}%
\DeclareMathOperator{\Coh}{Coh}%
\DeclareMathOperator{\gener}{gen}%

\newcommand{\nn}{\mathbb{N}}
\newcommand{\zz}{\mathbb{Z}}
\newcommand{\qq}{\mathbb{Q}}
\newcommand{\cc}{\mathbb{C}} 
\newcommand{\pp}{\mathbb{P}} 

\newcommand{\oo}{\mathcal{O}} 

\newcommand{\Xrig}{\mathcal{X}^{\mathrm{rig}}}
\newcommand{\Trig}{\mathcal{T}^{\mathrm{rig}}}
\newcommand{\frig}{f^{\mathrm{rig}}}

\newcommand{\DM}{Deligne-Mumford }
\newcommand{\bs}{\boldsymbol}
\newcommand{\aaa}{\mathbb{A}}

\def\cf{\textit{cf.}\kern.3em}
\def\eg{\textit{e.g.},\ }
\def\ie{\textit{i.e.},\ }
\def\resp{\textit{resp.}\kern.3em}
\renewcommand{\k}{\kern2pt}

\DeclareMathOperator{\red}{red}%
\DeclareMathOperator{\codim}{codim}%
\DeclareMathOperator{\rank}{rk}%
\DeclareMathOperator{\Hom}{Hom}%
\DeclareMathOperator{\Ext}{Ext}%
\DeclareMathOperator{\coker}{coker}%
\DeclareMathOperator{\diag}{diag}%
\DeclareMathOperator{\im}{Im}%
\DeclareMathOperator{\Id}{id}%
\DeclareMathOperator{\can}{can}%
\DeclareMathOperator{\Pic}{Pic}%
\DeclareMathOperator{\et}{\text{\'et}}%
\DeclareMathOperator{\Spec}{Spec}%
\DeclareMathOperator{\spec}{Spec}
\DeclareMathOperator{\Aut}{Aut}
\DeclareMathOperator{\rig}{rig}

\DeclareMathOperator{\tor}{tor}%
\DeclareMathOperator{\Mor}{Mor}%
\DeclareMathOperator{\ccs}{\cc^{\ast}}%
\renewcommand{\k}{\kern2pt}

\usepackage{amssymb,amsmath,amsfonts,amscd}

\begin{document}

  \title{Smooth toric Deligne-Mumford stacks }
 \author{Barbara Fantechi, Etienne Mann \and Fabio Nironi}
  \address{Barbara Fantechi, SISSA, Via Beirut 2-4, 34014 Trieste, Italy}
 \email{fantechi@sissa.it}
 \address{Etienne Mann, 
Universit\'e de Montpellier 2, D\'epartement de Math\'ematiques, Place Eug\`ene Bataillon,
F-34 095 Montpellier CEDEX 5, France}
 \email{etienne.mann@math.univ-montp2.fr}
  \address{Fabio Nironi, Department of Mathematics,
    Columbia University, 2990 Broadway MC 4418, New York, NY
    10027-6902, USA}
  \email{nironi@math.columbia.edu}
 \date{September 09 th 2009}
 
  \begin{abstract}
    We give a geometric definition of smooth toric \DM stacks using
    the action of a``torus''. We show that our definition is
    equivalent to the one of Borisov, Chen and Smith in terms of
    stacky fans. In particular, we give a geometric interpretation of
    the combinatorial data contained in a stacky fan.  We also give a
    bottom up classification in terms of simplicial toric varieties
    and fiber products of root stacks. 
  \end{abstract}

\maketitle



\tableofcontents
\newpage
\section*{Introduction}

\renewcommand\theequation{\Roman{equation}} A toric variety is a
normal, separated variety $X$ with an open embedding $T\hookrightarrow
X$ of a torus such that the action of the torus on itself extends to
an action on $X$. To a toric variety one can associate a fan, a
collection of cones in the lattice of one-parameter subgroups of $T$.
Toric varieties are very important in algebraic geometry, since
algebro-geometric properties of a toric variety translate in
combinatorial properties of the fan, allowing to test conjectures and
produce interesting examples.


In \cite{BCSocdms05} Borisov, Chen and Smith define toric \DM
stacks as explicit global quotient (smooth) stacks, associated to
combinatorial data called stacky fans. Later, Iwanari proposed in
\cite{Icts06} a definition of toric triple as an orbifold with a torus
action having a dense orbit isomorphic to the torus \footnotemark
\footnotetext{For the meaning of orbifold in this paper, see
  \S~\ref{subsec:smooth-dm-stacks}.} and he proved that the
2-category of toric triples is equivalent to the 2-category of
``toric stacks'' (We refer to \cite{Icts06} for the definition of
``toric stacks''). Nevertheless, it is clear that not all toric \DM stacks
are toric triples, since some of them are not orbifolds.




Then the generalization of the $\Delta$-collections defined for toric varieties
by Cox in \cite{Cox-functorofsmooth-1995} was done by Iwanari in
\cite{Iwanari-GeneralizationofCox-2007} in the orbifold case and by Perroni in
\cite{Perroni-notetoricDeligne-Mumford-2007} in the general case.

In this paper, we define a \DM torus $\TT$ as a Picard stack isomorphic to
$T\times\mathcal{B}G$, where $T$ is a torus, and $G$ is a finite abelian group; we then
define a smooth toric \DM stack as a smooth separated \DM stack with the
action of a \DM torus $\TT$ having an open dense orbit isomorphic to $\TT$.
We prove a classification theorem for
smooth toric \DM stacks and show that they coincide with those
defined by \cite{BCSocdms05}. 




The first main result of this paper is a bottom-up description of
smooth toric \DM stacks, as follows: the structure morphism $\XX\to X$
to the coarse moduli space factors canonically via the toric morphisms
\begin{displaymath}
\xymatrix{\XX \ar[r] & \Xrig \ar[r] &\XX^{\can}\ar[r] &X}
\end{displaymath}
where $\XX \to \Xrig$ is an abelian gerbe over $\Xrig$ 
; $\Xrig\to \XX^{\can}$
is a fibered product of roots of toric divisors
; and $\XX^{\can}\to X$ is the minimal orbifold having
$X$ as coarse moduli space. Here $X$ is a simplicial toric variety,
and $\Xrig$ and $\XX^{\can}$ are smooth toric \DM stacks.  More precisely, this
bottom up construction can be stated as follows.

\begin{thm}\label{thm:intro,charac}
Let $\XX$ be a smooth toric \DM stack with \DM torus isomorphic to
  $T\times \mathcal{B}G$. Denote by $X$ the coarse moduli space of
  $\XX$.  Denote by $n$ the number of rays of the fan of $X$.
  \begin{enumerate}
  \item There exist unique $(a_{1}, \ldots ,a_{n}) \in
    (\nn_{>0})^{n}$ such that the stack $\Xrig$ is isomorphic, as
    toric \DM stack, to
\begin{displaymath}
    \sqrt[a_{1}]{{D}^{\can}_{1}/\XX^{\can}}\times_{\XX^{\can}}\cdots\times_{\XX^{\can}}\sqrt[a_{n}]{{D}_{n}^{\can}/\XX^{\can}}.
  \end{displaymath} where
  $D_{i}^{\can}$ is the divisor corresponding to the ray $\rho_{i}$.
\item Given $G=\prod_{j=1}^{\ell} \mu_{b_{j}}$. There exist $L_{1},
  \ldots ,{L}_{\ell}$ in $\Pic(\Xrig)$ such that $\mathcal{X}$ is
  isomorphic, as toric \DM stack, to
  \begin{displaymath}
    \sqrt[b_{1}]{{L}_{1}/\Xrig}\times_{\Xrig}\cdots\times_{\Xrig}\sqrt[b_{\ell}]{{L}_{\ell}/\Xrig}.
  \end{displaymath} Moreover, for any $j\in\{1, \ldots ,\ell\}$, the
  class $[L_{j}]$ in $\Pic(\Xrig)/b_{j} \Pic(\Xrig)$
  is unique.
  \end{enumerate}
\end{thm}

In the process, we get a description of the Picard group of smooth
toric \DM stacks, which allows us to characterize weighted projective
stacks as complete toric orbifolds with cyclic Picard group (\cf
Proposition \ref{prop:charac,wps}). Moreover, we classify all complete
toric orbifolds of dimension $1$ (\cf Example \ref{exmp:toric,line}).
We also show that the natural map from the Brauer group of a smooth
toric \DM stack with trivial generic stabilizer to its open dense torus
is injective (\cf Theorem \ref{thm:gerbe,banale}).

The second main result of this article is to give an explicit relation
between the smooth toric \DM stacks and the stacky fans.
\begin{thm}\label{thm:intro,stacky}
  Let $\XX$ be a smooth toric \DM stack with coarse moduli space the
  toric variety denoted by $X$. Let $\Sigma$ be a fan of $X$ in
  $N_{\qq}:=N\otimes_{\zz} \qq$. Assume that the rays of $\Sigma$
  generate $N_{\qq}$. There exists a stacky fan such that $\XX$ is
  isomorphic, as toric \DM stack, to the smooth \DM stack associated
  to the stacky fan.  Moreover, if $\XX$ has a trivial generic
  stabilizer then the stacky fan is unique.
\end{thm}

When the smooth toric \DM stack $\XX$ has a generic stabilizer
the non-uniqueness of the stacky fan comes from three different
choices.  We refer to Remark \ref{rem:toric,choices} for a more
precise statement.
This result gives a geometrical interpretation of the combinatorial
data of the stacky fan.  In fact, the stacky fan can be read off the
geometry of the smooth toric \DM stack just like the fan can be read off the
geometry of the toric variety.
Notice that one can deduce the above theorem when $\XX$ is an orbifold
from Theorem 2.5 of \cite{Perroni-notetoricDeligne-Mumford-2007} and
Theorem 1.4 of \cite{Iwanari-GeneralizationofCox-2007} and the geometric 
characterization  of Theorem 1.3 in 
\cite{Iwanari-Logarithmic-2007} .

In the first part of this article, we fix the conventions and collect
some results on smooth \DM stacks, root constructions, rigidification,
toric varieties, Picard stacks and the action of a Picard stack.  In Section
\ref{sec:dm-tori}, we define \DM tori. Section
\ref{sec:defin-toric-orbif} contains the definition of smooth toric \DM
stacks. In Section \ref{sec:canon-toric-orbif}, we first define
canonical smooth \DM stacks and then we show that the canonical stack
associated to a simplicial toric variety is a smooth toric \DM stack (\cf
Theorem \ref{thm:charac,cano,orb}).  In Section
\ref{sec:reduc-toric-orbif}, we prove the first part of Theorem
\ref{thm:intro,charac}. In Section \ref{sec:toric-dm-stacks-1}, we
first prove in Proposition \ref{prop:ess,trivi,gerb} that the
essentially trivial banded gerbes over $\XX$ are in bijection with
finite extensions of the Picard group of $\XX$; then, we show that the
natural map from the Brauer group of a smooth toric \DM stack with trivial
generic stabilizer to its open dense torus is injective (\cf Theorem
\ref{thm:gerbe,banale}). Finally, we prove the second statement of
Theorem \ref{thm:intro,charac}. In Section
\ref{sec:toric-dm-stacks-2}, we prove Theorem \ref{thm:intro,stacky}
and give some explicit examples.  In Appendix
\ref{sec:action-picard-stack}, we have put some details about the
action of a Picard stacks.


\textbf{Acknowledgments} The authors would like to acknowledge support
from IHP, Mittag-Leffler Institut, SNS where part of this work was
carried out, as well as the European projects MISGAM and ENIGMA.  We
would like to thank Ettore Aldrovandi, Lev Borisov, Jean-Louis
Colliot-Th\'el\`ene, Andrew Kresch, Fabio Perroni, Ilya Tyomkin and
Angelo Vistoli for helpful discussions; in particular Aldrovandi for
explanations about group-stacks and reference \cite{Breen-Bitors-1990},
Borisov for pointing out a mistake in a preliminary version,
Colliot-Th\'el\`ene for \cite[\S 6]{Ggb68}, Tyomkin for \cite{Bal1993}
and Vistoli for useful information about the classification of gerbes.


\section{Notations and Background}
\numberwithin{equation}{section}

\subsection{Conventions and notations}\label{subsec:convention}

A scheme will be a separated scheme of finite type over $\cc$. A
variety will be a reduced, irreducible scheme. A point will be a
$\cc$-valued point.  The smooth locus of a variety $X$ will be denoted
by $X_{\sm}$.  

We work in the \'etale topology.  For an algebraic stack $\XX$, we
will write that $x$ is a point of $\XX$ or just $x\in \XX$ to mean
that $x$ is an object in $\XX(\cc)$; we denote by $\Aut(x)$ the
automorphism group of the point $x$.
We will say that a morphism between stacks is unique if it is unique up to a unique $2$-arrow.
As usual, we denote $\Gm$ the sheaf of invertible sections in
$\mathcal{O}_{\XX}$ on the \'etale site of $\XX$. 

\subsection{Smooth \DM stacks and orbifolds}\label{subsec:smooth-dm-stacks}

A \DM stack will be a separated \DM stack of finite type over $\cc$; we
will always assume that its coarse moduli space is a scheme.  An \textit{
  orbifold} will be a smooth \DM stack with trivial generic stabilizer.
For a smooth \DM stack $\XX$, we denote by $\eps_\XX$ or just $\eps$
the natural morphism from $\XX$ to its coarse moduli space $X$, which
is a variety with finite quotient singularities.



 Let $\iota:\UU\to\XX$ be an open embedding of
  irreducible smooth \DM stacks with complement of codimension at least
  $2$. We have that 
  \begin{itemize}
  \item The natural map $\iota^*:\Pic(\XX)\to\Pic (\UU)$ is an
  isomorphism. 
\item For any line bundle $L\in\Pic(\XX)$, the
  natural morphism $\iota^*:H^0(\XX,L)\to H^0(\UU,\iota^*L)$ is also an
  isomorphism.
\end{itemize}

The \textit{inertia stack}, denoted by $I(\XX)$, is defined to be the
fibered product $I(\XX):=\XX\times_{\XX\times\XX}\XX$. A point of
$I(\XX)$ is a pair $(x,g)$ with $x\in \XX$ and $g\in \Aut(x)$.  The
inertia stack of a smooth \DM stack is smooth but different components
will in general have different dimensions. The natural morphism
$I(\XX)\to\XX$ is representable, unramified, proper and a relative
group scheme. The identity section gives an irreducible component
canonically isomorphic to $\XX$; all other components are called \textit{
  twisted sectors}.  A smooth \DM stack of dimension $d$ is an orbifold
if and only if all the twisted sectors have dimension $\le d-1$, and
is \textit{canonical} if and only if all twisted sectors have dimension $\le
d-2$.

\begin{rem}[Sheaves on Global quotients]\label{rem:sheaves,global,quotient}
  According to \cite[Appendix]{Vitas-1989}, a coherent sheaf on a \DM stack
  $[Z/G]$ is a $G$-equivariant sheaf on $Z$ \ie the data of a coherent sheaf $L_{Z}$ on $Z$ and for every $g\in G$ an
  isomorphism $\varphi_{g}:L_{Z}\to g^{\ast} L_{Z}$ such that 
  $\varphi_{gh}=h^{\ast}\varphi_{g}\circ\varphi_{h}$.

  Notice that if $Z$ is a subvariety of $\cc^{n}$ of codimension higher or equal
  than two then an invertible sheaf on $[Z/G]$ is the structure sheaf
  $\mathcal{O}_{Z}$ and a one dimensional representation of $G$ \ie $\chi: G\to
  \cc^{\ast}$. A global section of such an invertible sheaf on $[Z/G]$ is a
  $\chi$-equivariant global section of $\mathcal{O}_{Z}$.
\end{rem}

We end this subsection with a proposition  extending to stacks a property of separated schemes. We will prove it in Appendix \ref{sec:appendix,proof}.

\begin{prop}\label{prop:sep,bis}
  Let $\XX$ and $\YY$ be two \DM stacks. Assume that $\XX$ is normal
  and $\YY$ is separated.  Let $\iota:\mathcal{U}\hookrightarrow \XX$ be a
  dominant open immersion of the \DM stack $\UU$. If $F,G :\XX \to \YY$ are
  two morphisms of stacks such that there exits a $2$-arrow $F\circ
  \iota \stackrel{\beta}{\Rightarrow}G\circ \iota$ then there exists a
  unique $2$-arrow $\alpha:F\Rightarrow G$ such that $\alpha\ast
  \Id_{\iota}=\beta$.
\end{prop}

The previous proposition is well-known for $\XX$ a reduced scheme and
$\YY$ a separated scheme. Nevertheless, if $\XX$ is not a normal stack
we have the following counter-example:
Let $\YY$ be $\mathcal{B}\mu_{2}$. Let $\XX$ be a rational curve with one node. Let
$F_{1}:\XX\to \YY$ (resp. $F_{2}$) be a stack morphism given by a non
trivial (resp. trivial) double cover of $\XX$. Putting
$\UU=\XX\setminus \{\textrm{node}\}$, the proposition is false. 

\subsection{Root constructions}\label{subsec:root-construction}
For this subsection we refer to the paper of Cadman
\cite{Cstc-2007} (see also  \cite[Appendix B]{AGVgwdms}).
In this part $\XX$ will be a \DM stack over $\cc$ (it is enough to
assume that $\XX$ is Artin.)

\subsubsection{Root of an invertible sheaf}\label{sec:root-line}
This part follows closely Appendix B of \cite{AGVgwdms}.
Let $L$ be an invertible sheaf on the \DM stack $\XX$. Let $b$ be a
positive integer.  We denote by $\sqrt[b]{L/\XX}$ the following fiber
product
 \begin{displaymath}
   \xymatrix{\sqrt[b]{L/\XX} \ar[r] \ar[d] \ar@{}[dr]|{\square}& \mathcal{B}\ccs
   \ar[d]^{\wedge b} \\ \XX \ar[r]^-{L}&\mathcal{B}\ccs }
 \end{displaymath}
 where $\wedge b :\mathcal{B}\ccs \to \mathcal{B}\ccs$ sends an
 invertible sheaf $M$ over a scheme $S$ to $M^{\otimes b}$. More explicitly,
 an object of $\sqrt[b]{L/\XX}$ over $f:S\to \XX$ is a couple
 $(M,\varphi)$ where $M$ is an invertible sheaf $M$ on the scheme $S$ and
 $\varphi:M^{\otimes b}\stackrel{\sim}{\to} f^{\ast}L$ is an isomorphism. The arrows are defined
 in an obvious way. 
 
 The morphism $\sqrt[b]{L/\XX}\to B\ccs$ corresponds to an invertible
 sheaf, denoted by $L^{1/b}$ in \cite{BCqh-2007}, on $\sqrt[b]{L/\XX}$
 whose $b$-th power is isomorphic to the pullback of $L$.

The stack $\sqrt[b]{L/\XX}$ is a $\mu_{b}$-banded
 gerbe over $\XX$ (see the second paragraph of Subsection \ref{subsec:gerbe+root} below).  The Kummer exact sequence
\begin{displaymath}
  \xymatrix{1\ar[r]&\mu_{b}\ar[r]&\Gm\ar[r]^{\wedge b}&\Gm\ar[r]&1}
\end{displaymath} 
induces the boundary morphism $\partial :H^{1}_{\et}(\XX,\Gm)
\to H^{2}_{\et}(\XX,\mu_{b})$.
The cohomology class of the $\mu_{b}$-banded gerbe $\sqrt[b]{L/\XX}$ in
$H^{2}_{\et}(\XX,\mu_{b})$ is the image by $\partial$ of the class $[L]\in
H^{1}_{\et}(\XX,\Gm)$. 

The gerbe is trivial if and only if
the invertible sheaf $L$ has a $b$-th root in $\Pic(\XX)$.
More generally, the gerbe $\sqrt[b]{L/\XX}$ is isomorphic, as a
$\mu_{b}$-banded gerbe, to $\sqrt[b]{L'/\XX}$
if and only if $[L]=[L']$ in $\Pic(\XX)/b \Pic(\XX)$.
We have the following morphism of short exact sequences:
 \begin{equation}\label{eq:21}
     \xymatrix{0\ar[r]& \zz\ar[r]^{\times b}\ar[d]& \zz\ar[r]\ar[d]&
       \zz/b\zz \ar[r] \ar@{=}[d]&0\\
 0\ar[r]& \Pic(\XX)\ar[r]& \Pic(\sqrt[b]{L/\XX})\ar[r]& \zz/b\zz \ar[r]&0}
   \end{equation}
   where the first and second vertical morphisms are defined by $1\mapsto {L}$
   and $1\mapsto L^{1/b}$, respectively.
\subsubsection{Roots of  effective Cartier divisors}\label{subsubsec:root-line-section}
In the articles \cite{Cstc-2007} and \cite{AGVgwdms}, the authors
define the notion of root of an invertible sheaf with a section on an
algebraic stack: here, we only consider roots of effective Cartier
divisors on a smooth algebraic stack, since this is what we will use.

Let $n$ be a positive integer. Consider the quotient stack
$[\aaa^{n}/(\ccs)^{n}]$ where the action of $(\ccs)^{n}$ is given
multiplication coordinates by coordinates. Notice that
$[\aaa^{n}/(\ccs)^{n}]$ is the moduli stack of $n$ line bundles with
$n$ global sections.  Let   $\bs{a}:=(a_{1}, \ldots
,a_{n})\in (\nn_{>0})^{n}$ be a $n$-tuple. Denote by $\wedge
\bs{a}: [\aaa^{n}/(\ccs)^{n}] \to [\aaa^{n}/(\ccs)^{n}]$ the stack morphism
defined by sending $x_{i}\mapsto x_{i}^{a_{i}}$ and
$\lambda_{i}\mapsto \lambda_{i}^{a_{i}}$ where $x_{i}$ (resp. $\lambda_{i}$)
are coordinates  of $\aaa^{n}$
(resp. $(\ccs)^{n}$).

 Let $\XX$ be a smooth algebraic stack. 
 Let $\bs{D}:=(D_{1}, \ldots ,D_{n})$ be $n$
effective Cartier divisors. 
  The $\bs{a}$-th root of
$(\XX,\bs{D})$ is the fiber product
\begin{displaymath}
  \xymatrix{\sqrt[\bs{a}]{\bs{D}/\XX} \ar@{}[rd]|{\square}\ar[r]\ar[d]_{\pi}& [\aaa^{n}/(\ccs)^{n}]
  \ar[d]^{\wedge \bs{a}} \\ \XX \ar[r]^{\bs{D}}& [\aaa^{n}/(\ccs)^{n}]}
\end{displaymath}
The morphism $\sqrt[\bs{a}]{\bs{D}/\XX}\to  [\aaa^{n}/(\ccs)^{n}]$
corresponds to the effective Cartier divisors
$\widetilde{\bs{D}}:=(\widetilde{D}_{1}, \ldots ,\widetilde{D}_{n})$,
where $\widetilde{D}_{i}$ is the reduced closed substacks
$\pi^{-1}(D_{i})_{\red}$. More explicitly, an object of
$\sqrt[\bs{a}]{\bs{D}/\XX}$ over a scheme $S$ is a couple
$(f,(\widetilde{D}_{1}, \ldots ,\widetilde{D}_{n}))$ where $f:S \to
\XX$ is a morphism and for any $i$, $D_{i}$ is an effective divisor on $S$
such that $a_{i}\widetilde{D}_{i}=f^{\ast}D_{i}$.

We have the following properties :
\begin{enumerate}
\item\label{item:11} The fiber product of $\sqrt[a_{i}]{D_{i}/\XX}$
  over $\XX$ is isomorphic to $\sqrt[\bs{a}]{\bs{D}/\XX}$ (\cf Remark
  2.2.5 of \cite{Cstc-2007}).
\item\label{item:12}  The canonical morphism $\sqrt[\bs{a}]{\bs{D}/\XX} \to \XX$ is an
  isomorphism over $\XX\setminus\cup_{i}D_{i}$.
\item\label{item:9} If $\XX$ is smooth, each $D_{i}$ is smooth and  the $D_{i}$ have simple
  normal crossing then $\sqrt[\bs{a}]{\bs{D}/\XX}$ is smooth (\cf
  Section 2.1 of \cite{BCqh-2007}) and $\widetilde{D}_{i}$ have simple
  normal crossing.
\item\label{item:13} We have the following morphism of short exact
  sequences (\cf Corollary 3.1.2 \cite{Cstc-2007})
\begin{equation}\label{eq:29}
     \xymatrix@1{0\ar[r]&\zz^{n}\ar[r]^{\times \bs{a}} \ar[d]&
       \zz^{n}\ar[r] \ar[d]& \prod_{i=1}^{n}\zz/a_{i}\zz\ar[r]\ar@{=}[d]&0
 \\0\ar[r]&\Pic(\XX)\ar[r]^(.4){\pi^*} &
\Pic\left(\sqrt[\bs{a}]{\bs{D}/\XX}\right)\ar[r]^{q} 
& \prod_{i=1}^{n}\zz/a_{i}\zz\ar[r]&0}
   \end{equation}
 where the first and second vertical morphisms are defined by
 $e_{i}\mapsto \OO(D_{i})$ and  $e_{i}\mapsto
 \OO(\widetilde{D}_{i})$, respectively.
   Every invertible sheaf $L\in
   \Pic\left(\sqrt[\bs{a}]{\bs{D}/\XX}\right) $ can be written in a
   unique way as $L\cong\pi^*
   M\otimes\prod_{i=1}^n\mathcal{O}(k_i\widetilde{D}_i)$ where
   $M\in\Pic(\XX)$ and $0\leq k_i<a_{i}$; the morphism $q$ maps $L$
   to $(k_1,\ldots,k_n)$.

\end{enumerate}

We finish this section with the following observation.
 Let $D_{1}$ and $D_{2}$ be two effective Cartier divisors on
  $\XX$ such that $D_{1}\cap D_{2} \neq \emptyset$. The stacks
  $\sqrt[a]{D_{1}\cup D_{2}/\XX}$ and
  $\sqrt[(a,a)]{(D_{1},D_{2})/\XX}$ are not isomorphic.  Indeed, the
  stabilizer group at any point in the preimage of $x\in D_{1}\cap
  D_{2}$ in $\sqrt[a]{D_{1}\cup D_{2}/\XX}$ (resp.
  $\sqrt[(a,a)]{(D_{1},D_{2})/\XX}$) is $\mu_{a}$ (resp.
  $\mu_{a}\times \mu_{a}$).

\subsection{Rigidification}\label{subsec:rigidification}
 \def\IgenX{I^{\gener}(\XX)} \def\IgenT{I^{\gener}(\TT)}
 
 In this subsection, we sum up some results on the rigidification of an
 irreducible $d$-dimensional smooth \DM stack $\XX$. Intuitively, the
 rigidification of $\XX$ by a central subgroup $G$ of the
 \textit{generic stabilizer}  is constructed as
 follows: first, one constructs a prestack where the objects are the
 same and the automorphism groups of each object $x$ are the quotient
 $\Aut_{\XX}(x)/G$; then the rigidification $\XX \fs G$ is the
 stackification of this prestack. For the most general construction we
 refer to \cite[Appendix A]{abramovich-2007} (see also \cite[Section
 5.1]{ACVtbac-2003}, \cite{Rgasa-2005} and \cite[Appendix
 C]{AGVgwdms}).
 
We  consider the union $\IgenX\subset I(\XX)$
 of all $d$-dimensional components of $I(\XX)$; it is a subsheaf of
 groups of $I(\XX)$ over $\XX$ which is called the generic stabilizer.
 Most of the time in this article, we will rigidify by the 
 generic stabilizer. In this case, we write $ \Xrig$ in order to mean
 $\XX \fs \IgenX$ and call it \textbf{the} rigidification.

The rigidification $r:\XX\to \Xrig$ has the following properties: 
\begin{enumerate}
\item the coarse moduli space of $\Xrig$ is the coarse moduli space of $\XX$,
\item $\Xrig$ is an orbifold,
\item if $\XX$ is an orbifold then $\Xrig$ is $\XX$,
\item the morphism $r$ makes $\XX$ into a gerbe over $\Xrig$.
\end{enumerate}

We refer to Theorem 5.1.5.(2) of \cite{ACVtbac-2003} for the proof of
the following proposition.
\begin{prop}[Universal property of the rigidification]\label{prop;universal}
  Let $\XX$ be a smooth \DM stack. Let $\YY$ be an orbifold. Let
  $f:\XX\to \YY$ be a dominant stack morphism. Then there exists
  $g:\Xrig \to \YY$ and a $2$-morphism $\alpha:g\circ r \Rightarrow f$
  such that the following is $2$-commutative
  \begin{displaymath}
    \xymatrix{\XX \ar[r]^{r} \ar[dr]_{f} & \Xrig\ar[d]^{\exists \k g }
    \\ & \YY}\end{displaymath}
If there exists $g':\Xrig \to
  \YY$ and a $2$-morphism $\alpha':g'\circ r \Rightarrow f$ satisfying
    the same property then there exists a unique $\gamma:g'\Rightarrow
    g$ such that $\alpha\circ(\gamma\ast \Id_{r})=\alpha'$.
\end{prop}

\subsection{Diagonalizable group schemes}\label{subsec:diag-group-scheme}

In this short subsection, we recall some results on
diagonalizable groups.

\begin{defn} A group scheme $G$ over $\spec \cc$ will be called
  \textit{diagonalizable} if it is isomorphic to the product of a torus
  and a finite abelian group.
\end{defn}

We use multiplicative notation for diagonalizable group.  For any
diagonalizable group $G$, its character group $G^\vee:=\Hom(G,\cc^{\ast})$ is
a finitely generated abelian group (or coherent $\zz$-module).  The
duality contravariant functor $G\mapsto G^\vee$ induces an equivalence
of categories from diagonalizable to coherent $\zz$-module. Its
inverse functor is given by $F\mapsto G_F:=\Hom(F,\cc^{\ast})$. Both
$G\mapsto G^\vee$ and $F\mapsto G_F$ are contravariant and exact.

\subsection{Toric varieties}\label{subsec:geometry-toric}
We recall some results on toric varieties that can be found in
\cite{Fitv} (see also \cite{Cltv2000}). The principal construction
 used in this paper is the description of toric varieties as
global quotients found by Cox (see \cite{Chcr95}).  

We fix a torus $T$, and denote by $M=T^{\vee}$ the lattice of
characters and by $N:=\Hom(M,\zz)$ the lattice of one-parameter
subgroups.  A toric variety $X$ with torus $T$ corresponds to a fan
$\Sigma(X)$, or just $\Sigma$, in $N_\qq:=\nn\otimes_{\zz}\qq$, which
we will always assume to be simplicial.

Let $\rho_{1}, \ldots ,\rho_{n}$ be the one-dimensional cones, called
rays, of $\Sigma$. For any ray $\rho_{i}$, denote by $v_{i}$ the unique
generator of $\rho_{i}\cap N$.  For any $i$ in $\{1, \ldots ,n\}$, we
denote by $D_{i}$ the irreducible $T$-invariant Weil divisor defined by
the ray $\rho_{i}$. The free abelian group of $T$-invariant Weil
divisor is denoted by $L$. 

Let $\iota:M\to L$ be the morphism that sends $m\mapsto \sum_{i=1}^{n}
m(v_{i})$. If the rays span $\nn_{\qq}$ (which is not a strong
assumption \footnote{\label{footnote} Indeed, if the rays do not span $N_{\qq}$ then $X$ is
  isomorphic to the product of a torus and a toric variety
  $\widetilde{X}$ whose rays span $\widetilde{N}_{\qq}$.}), the morphism $\iota$ is injective, and fits into an exact
sequence in $\Coh(\zz)$
\begin{align}\label{eq:9}
\xymatrix{ 0\ar[r]& M\ar[r]^-{\iota} & L \ar[r]& 
A \ar[r]& 0.}
\end{align} where $A$ is the class group of $X$ (\ie the
Chow group $A^{1}(X)$).
We deduce that the short exact sequence of diagonalizable groups
\begin{align}\label{eq:10}
\xymatrix{1\ar[r]& G_{A}\ar[r]& G_{L}\ar[r]& T\ar[r]& 1.}
\end{align}

Let $Z_{\Sigma}\subset \cc^n$ be the $G_{L}=(\cc^{\ast})^n$-invariant
open subset defined as $Z_{\Sigma}:=\cup_{\sigma\in \Sigma} Z_\sigma$,
where $Z_\sigma:=\{x\,|\,x_i\ne 0\text{ if }\rho_{i}\notin\sigma\}$.
The induced action of $G_A$ on $Z_{\Sigma}$ has finite stabilizers (by
the simpliciality assumption) and $X$ is the geometric quotient
$Z_{\Sigma}/G_{A}$, with torus $(\cc^{\ast})^{n}/G_{A}$ (see Theorem
2.1 of \cite{Chcr95}).  For any $i\in\{1, \ldots ,n\}$, the
$T$-invariant Weil divisor $D_i\subset X$ is the geometric quotient
\begin{equation}\label{eq:15}
  \left(\{x_{i}=0\}\cap Z_{\Sigma}\right)/G_{A}. 
\end{equation}

If $X$ is smooth then the natural morphism $L\to\Pic(X)$
given by $e_i\mapsto \oo_X(D_i)$ is surjective and has kernel
$M$; in other words, it induces a natural isomorphism $A\to \Pic(X)$.

 If $X$ is a $d$-dimensional toric variety, we will write $X^0$ for the
 union of the orbits of dimension $\ge d-1$; in other words, $X^0$ is
 the toric variety associated to the fan $\Sigma_{\le 1}:=\{\sigma\in
 \Sigma\,| \dim\sigma\le 1\}$. The toric variety $X^0$ is always smooth
 and the toric divisors $D_\rho^0$ are smooth, disjoint, and homogeneous
 under the $T$-action (with stabilizer the one-dimensional subgroup
 which is the image of $\rho$).

\subsection{Picard stacks and action of a Picard stack}
Deligne defined Picard stacks in \cite[Expos\'e XVIII]{SGA4} as stacks
analogous to sheaves of abelian groups. For the reader's convenience, we
collect here a sketch of the definition and the main properties we
need; details can be found in \cite[Expos\'e XVIII]{SGA4} and also in
\cite[Section 14]{LMBca}.

Here we summarize the definition of a Picard stack. For the details we
refer to Definition \ref{defi:picard,stacks}. 
\begin{defn} Let $\GG$ be a stack over a base scheme $S$. 
  A \textit{Picard stack} $\GG$ over  $S$ is given by the
  following set of data:
  \begin{itemize}
  \item a multiplication stack morphism $m:\GG\times\GG\to\GG$, also
    denoted by $m(g_1,g_2)=g_1\cdot g_2$;
  \item an associativity $2$-arrow $(g_1\cdot g_2)\cdot g_3\Rightarrow
    g_1\cdot(g_2\cdot g_3)$;
  \item a commutativity $2$-arrow $g_1\cdot g_2\Rightarrow g_2\cdot
    g_1$.
  \end{itemize}
These data satisfy some compatibility relations, which we
list in  \ref{defi:picard,stacks}.
\end{defn}

The definition implies that there also exists an identity $e:S\to\GG$
and an inverse $i:\GG\to\GG$ with the obvious properties; in
particular, a $2$-arrow $\eps:(e\cdot g)\Rightarrow g$.

\begin{defn}[See Section 1.4.6 in \cite{SGA4}]  Let $\GG$, $\GG '$ be two Picard stacks. A \textit{morphism of Picard stacks} $F:\GG\to \GG'$ is a morphism of stacks
  and a $2$-arrow $\alpha$ such that for any two objects $g_{1},g_{2}$
  in $\GG$, we have
\begin{displaymath}
      F(g_{1}\cdot g_{2}) \stackrel{\alpha}{\Rightarrow }F(g_{1})\cdot
      F(g_{2}).   
 \end{displaymath}
\end{defn}

Again we refer to  Appendix B for the list of compatibilities satisfied by $\alpha$.
The Picard stacks over $S$ form a category where the objects are
Picard stacks and morphisms are equivalence classes of morphism of
Picard stacks.

\begin{rem}\label{rem:cpx,picard} 
  To any complex $G^\bullet:=[G^{-1}\to G^0]$ of sheaves of abelian
  groups, we can associate a Picard stack $\GG$.  In this paper,
  $G^{\bullet}$ will be a complex of diagonalizable groups and the
  associated Picard stack is the quotient stack $[G^{-1}/G^{0}]$.
\end{rem}

\begin{prop}[See Proposition 1.4.15 in
  \cite{SGA4}]\label{prop:Picard,stack} The functor that associates to
  a length $1$ complex of sheaves of abelian groups a Picard stack
  induces an equivalence of categories between the derived category,
  denoted by $D^{[-1,0]}(S,\zz)$, of
  length $1$ complexes of sheaves of abelian groups and the category
  of Picard stacks.
\end{prop}

In particular, if $G$ is any sheaf of abelian groups on the base
scheme $S$, the quotient
$[S/G]$, i.e. the gerbe $\mathcal{B}G$, is naturally a Picard stack.

We finish this section with a sketch of the definition of an action of
a Picard stack on a stack. This is a generalization of the action of a
group scheme on a stack defined by Romagny in \cite{Rgasa-2005}. We
refer to Definition \ref{defi:action,picard,stack} for the details.

\begin{defi}[Action of a Picard stack]\label{defi:action,picard,stack,intro}
  Let $\mathcal{G}$ be a Picard stack. Denote by $e$
  the neutral section and by $\epsilon$ the corresponding $2$-arrow.
  Let $\XX$ be a stack. An \textit{action} of $\mathcal{G}$ on
  $\mathcal{X}$ is the following data :
  \begin{itemize}
  \item  a stack morphism $a : \GG \times \XX \to \XX$, also denoted
  by $a(g,x)=g\times x$;
\item  a 2-arrow $e\times x \Rightarrow x$;
\item  an associativity  2-arrow $(g_{1}\cdot g_{2})\times x\Rightarrow g_{1}\times(g_{2}\times x)$.
\end{itemize}
These data satisfy some compatibility relations, which we list in Appendix B.
\end{defi}

\section{\DM tori}\label{sec:dm-tori}
In this Section we define \DM tori which will play the role of the
torus for a toric variety.

We start with a technical Lemma.
\begin{lem}\label{lem:dervied,cat} 
  Let $\phi:A^0\to A^1$ be a morphism of finitely generated abelian groups such
  that $\ker \phi$ is free. In the
  derived category of complexes of finitely generated abelian groups
  of length $1$, the complex $[A^0\to A^{1}]$ is isomorphic to $[\ker
  \phi\stackrel{0}{\to}\coker\phi]$.
\end{lem}

\begin{proof}
  We have a  morphism of complexes
 $$[A^0\stackrel{\phi}{\to} A^1] \to
  [A^0/A^{0}_{\tor}\stackrel{\widetilde{\phi}}{\to}
  A^1/A^{0}_{\tor}]$$
 induced by the quotient morphisms. As $\ker \phi$
  is free, we deduce after a diagram chasing that this morphism is a quasi-isomorphism of
  complexes.  In  the derived category, we replace
  $A^1/A^{0}_{\tor}$ with a projective resolution
  $[\zz^{\ell}\stackrel{Q}{\to}\zz^{d+\ell}]$. Then the mapping cone of the
  morphism of complexes $[0\to A^0/A^{0}_{\tor}] \to [Q:\zz^{\ell}\to \zz^{d+\ell}]$
is $\left[[\widetilde{\phi}Q]:A^0/A^{0}_{\tor}\times
  \zz^{\ell}\to\zz^{d+\ell}\right]$ which is quasi-isomorphic to
$[A^0/A^{0}_{\tor}\stackrel{\widetilde{\phi}}{\to} A^1/A^{0}_{\tor}]$.
A morphism of free abelian groups $f$ is quasi-isomorphic to the
complex $[\ker f \stackrel{0}{\to} \coker f]$ and this finishes the proof. 
\end{proof}

The reader who is familiar with the article \cite{BCSocdms05} has
 probably recognized part of the construction of the stack associated to
a stacky fan.

\begin{rem}
  Let $\phi:A^0\to A^1$ be a morphism of finitely generated abelian
  groups as in the above lemma. Applying the contravariant functor
  $\Hom(\cdot,\cc^{\ast})$ of Section \ref{subsec:diag-group-scheme}
  to the complex $A^{0}\to A^{1}$, we get a length $1$ complex of
  diagonalizable groups $[G_{A^{1}}\xrightarrow{G_\phi} G_{A^{0}}]$. According to
  Remark \ref{rem:cpx,picard}, the associated Picard stack
  $[G_{A^{0}}/G_{A^{1}}]$ is a \DM stack if and only if the cokernel of
  $\phi$ is finite.
\end{rem}

\begin{exmp}\label{exmp:P(w),torus} Let $w_{0}, \ldots ,w_{n}$ be in
  $\nn_{>0}$. Let $\phi:\zz^{n+1}\to \zz$ that sends $(a_{0}, \ldots ,a_{n})$ to
  $\sum w_{i}a_{i}$. We have that $\ker\phi=\zz^{n}$ and $\coker\phi=\zz/d\zz$
  where $d:=\gcd(w_{0}, \ldots , \ldots ,w_{n})$. Hence, the associated Picard
  stack is $(\cc^{\ast})^{n}\times \mathcal{B}\mu_{d}$.
\end{exmp}

\begin{defn} A \textit{\DM torus} is a Picard stack over $\spec \cc$ which
  is obtained as a quotient $[G_{A^{0}}/G_{A^{1}}]$, where
  $\phi:A^0\to A^1$ is  a morphism of finitely generated abelian groups such
 that $\ker \phi$ is free and $\coker \phi$ is finite.
\end{defn}

Let $G$ be a finite abelian group. Notice that $\mathcal{B}G$ is a \DM
torus. Recall that by Proposition
\ref{prop:Picard,stack}, $T\times \mathcal{B}G$ has a natural structure
of Picard stack.
\begin{defn}
  A short exact sequence of Picard $S$-stacks is the sequence of morphisms of Picard $S$-stacks associated to a distinguished triangle in $D^{[-1,0]}(S)$. 
\end{defn}
\begin{prop}\label{prop:T,BG}
  Any \DM torus $\TT$ is isomorphic as Picard stack to
  $T\times\mathcal{B}G$ where $T$ is a torus and $G$ is a finite
  abelian group.
\end{prop}

\begin{proof}
Let $\TT=[G_{A^0}/G_{A^1}]$ with $\phi:A^0\to A^1$ as above. The distinguished triangle $[\ker G_\phi\to 0]\to [G_{A^1}\xrightarrow{G_\phi} G_{A^0}]\to [0\to\coker G_{\phi}]$ in the derived category $D^{[-1,0]}(\spec\mathbb{C})$ induces an exact sequence of Picard stacks $1\to \mathcal{B}G\to\TT\to T\to 1$ where $T:=G_{A^0}/G_{A^1}$. Proposition \ref{prop:Picard,stack} and Lemma \ref{lem:dervied,cat} imply that there is a non canonical isomorphism of Picard stacks $\TT\equiv \mathcal{B}G\times T$. 
\end{proof}
Note that the scheme $T$ in the previous proof is the coarse moduli space of $\TT$.



\section{Definition of toric \DM stacks}\label{sec:defin-toric-orbif}

\begin{defi}\label{defi:toric,orbifold}
\textit{A smooth toric \DM stack} is a smooth separated \DM stack
$\mathcal{X}$ together with an open immersion of a \DM torus
$\iota:\mathcal{T}\hookrightarrow \mathcal{X}$ with dense image such that
the action of $\mathcal{T}$ on itself extends to an action
$a:\mathcal{T}\times\mathcal{X}\to\mathcal{X}$.
\end{defi}

As in this paper all toric \DM stacks are smooth, we will write toric \DM
stack instead of smooth toric \DM stack. We will see later in Theorem \ref{thm:toric,BCS} that our definition a posteriori coincide with that in \cite{BCSocdms05} via stacky fans.
It seems natural to define a toric \DM stack by replacing smooth with normal in the above definition. All the definitions and results in this section apply also in this case, with the exception of Proposition \ref{prop:toric,orbifold,cms,toric,variety} and Lemma \ref{lem:igen,G,per,X}. Ilya Tyomkin is currently working on this. 
\textit{A toric orbifold} is a toric \DM stack with generically
trivial stabilizer. A toric \DM stack is a toric orbifold if and only if its \DM torus is an ordinary torus. Hence, the notion of toric orbifold is the same as the one
used in Theorem 1.3 of \cite{Icts06}.

\begin{rem}
\begin{enumerate}
\item Separatedness of $\XX$ and Proposition \ref{prop:sep,bis} imply that
  the action of $\TT$ on $\XX$ is uniquely determined by $\iota$.
\item Notice that we have assumed in Section \ref{subsec:smooth-dm-stacks} that the coarse
  moduli space is a scheme. Without this assumption, if the coarse
  moduli space $X$ of a toric \DM stack is a smooth and complete algebraic
  space then the main theorem of Bialynicki-Birula in \cite{Bal1993}
  implies that $X$ is a scheme. We don't know whether such an assumption is necessary in general.
\item A toric variety admits a structure of toric \DM stack if and only
  it is smooth.
\end{enumerate}
\end{rem}

\begin{prop}\label{prop:action,morphism}
  Let $\XX$ be a smooth \DM stack  together with an open dense
  immersion of a \DM torus $\iota:\mathcal{T}\hookrightarrow
  \mathcal{X}$ such that the action of $\mathcal{T}$
  on itself extends to a stack morphism
  $a:\mathcal{T}\times\mathcal{X}\to\mathcal{X}$. Then the stack
  morphism $a$ induces naturally an action of $\TT$ on $\XX$.
\end{prop}

\begin{proof}  We will define a
  $2$-arrow $\eta : a\circ (e,\Id_{\XX})\Rightarrow \Id_{\XX}$ and a
  $2$-arrow $\sigma: a \circ (m,\Id_{\XX}) \Rightarrow a \circ
  (\Id_{\XX},a)$ such that they verify Conditions
  (\ref{eq:pentagon.bis}) and (\ref{eq:compatibility.neutral.action})
  of Definition \ref{defi:action,picard,stack}. 
  We will only  prove the existence of $\eta$
  because the existence of $\sigma$ and  the relations
  (\ref{eq:pentagon.bis}) and (\ref{eq:compatibility.neutral.action})
  follow with a similar argument.

  Denote by $e:\spec \cc \to \TT$ the neutral element of $\TT$ and by
  $m:\TT\times \TT \to \TT$ the multiplication on $\TT$.  Denote by
  $\varepsilon$ the $2$-arrow $ m\circ (e,\Id_{\TT})\Rightarrow
  \Id_{\TT}$. As the stack morphism $a$ extends $m$, we have a
  2-arrow $\alpha:a\circ(\Id_{\TT},\iota) \Rightarrow \iota\circ m$.
  Denote by $\beta$ the 2-arrow $(e,\Id_{\XX})\circ \iota\Rightarrow
  (\Id_{\TT},\iota)\circ (e,\Id_{\TT})$. Consider the two stack
  morphisms :
  \begin{displaymath}
    \xymatrix{\TT \ar@^{^{(}->}[r]^{\iota}& \XX\ar@/_1pc/[r]_{a\circ(e,\Id_{\XX})} \ar@/^1pc/[r]^{\Id_{\XX}}& \XX.}
  \end{displaymath}
  Applying Proposition \ref{prop:sep,bis} with the composition of the
  following 2-arrows
\begin{displaymath}
  \xymatrix{
a\circ(e,\Id_{\XX})\circ\iota \ar@2{->}[r]^-{\Id_{a}\ast\beta}& a\circ
  (\Id_{\TT},\iota)\circ(e,\Id_{\TT})\ar@2{->}[r]^-{\alpha\ast\Id_{(e,\Id_{\TT})}}
& \iota\circ m\circ
  (e,\Id_{\TT})\ar@2{->}[r]^{\Id_{\iota}\ast\varepsilon}
&\iota\circ\Id_{\TT}=\Id_{\XX}\circ\iota,}
\end{displaymath}
we deduce the existence of $\eta: a\circ (e,\Id_{\XX})\Rightarrow
\Id_{\XX}$.
\end{proof}

\begin{defn} Let $\XX$ (resp. $\XX'$) be a toric \DM stack with \DM
  torus $\mathcal{T}$ (resp. $\mathcal{T}'$).  A \textit{morphism of
    toric \DM stacks} $F:\XX \to \XX'$ is a morphism of stacks between
  $\XX$ and $\XX'$ which extends a morphism of \DM tori $\TT\to\TT'$:
\end{defn}
\begin{rem} The extended morphism $F$ in the previous definition is unique by Proposition \ref{prop:sep,bis}. Moreover the definition of morphism between Picard stacks and 
  Proposition \ref{prop:sep,bis} provide us  the following $2$-cartesian diagram:
  \begin{displaymath}
    \xymatrix{\XX\times \TT \ar[rr]^-{(F,F\vert_{\TT})} \ar[d]_-{a} \ar@{}[rrd]|{\square}&&
    \XX'\times \TT' \ar[d]^-{a'} \\ \XX\ar[rr]^{F} && \XX'}
  \end{displaymath}
\end{rem}

\begin{prop}\label{prop:toric,orbifold,cms,toric,variety} 
  Let $\XX$ be a toric \DM stack with \DM torus $\mathcal{T}$. Let $X$
  (resp. $T$) be the coarse moduli space of $\mathcal{X}$ (resp.
  $\TT$).  Then $X$ has a structure of simplicial toric variety with
  torus $T$ where the open dense immersion
  $\overline{\iota}:{T}\hookrightarrow {X}$ and the action
  $\overline{a}:T\times X \to X$ is induced respectively by
  $\iota:\mathcal{T}\hookrightarrow \mathcal{X}$ and
  $a:\mathcal{T}\times \mathcal{X}\to \XX$.
\end{prop}

\begin{proof} The morphisms $\iota$
  and $a$ induce morphisms on the coarse moduli spaces $\bar\iota:T\to
  X$ and $\bar a:T\times X\to X$, by the universal property of the
  coarse moduli space. It is immediate to verify that $\bar\iota$ is
  an open embedding with dense image and $\bar a$ is an action,
  extending the action of $T$ on itself. On the other hand, since $X$
  is the coarse moduli space of $\XX$, it is a normal separated
  variety with finite quotient singularities. Therefore $X$ is a toric
  variety, and it is simplicial by \cite[\S 7.6 p.121]{Hms2003}(see
  also \cite[Theorem 3.1 p.28]{Cltv2000}).
\end{proof}

\begin{rem}[Divisor multiplicities]\label{rem:mult}
  According to \cite[Corollary 5.6.1]{LMBca}, the structure morphism
  $\eps:\XX\to X$ induces a bijection on reduced closed substacks.
  For each $i=1,\ldots,n$, denote by $\DD_i\subset \XX$ the reduced
  closed substack with support $\varepsilon^{-1}(D_{i})$.  Since
  $D_i\cap X_{\sm}$ is a Cartier divisor, there exists a unique
  positive integer $a_i$ such that $\eps^{-1}(D_i\cap
  X_{\sm})=a_i(\DD_i\cap \eps^{-1}(X_{\sm}))$.  We call
  $a=(a_1,\ldots,a_n)$ \textit{the divisor multiplicities of $\XX$}.
\end{rem}

Let $\XX$ be a  toric \DM stack with \DM torus $\TT=T\times
\mathcal{B}G$. By Appendix \ref{sec:action-picard-stack}, we have that
$\mathcal{B}G$ acts on $\XX$. Proposition \ref{prop:action,morphism,BG}
implies that we have an \'etale morphism $j:G\times \XX\to \IgenX$.

\begin{lem}\label{lem:igen,G,per,X}
  Let $\XX$ be a  toric \DM stack with \DM torus
  $\TT=T\times \mathcal{B}G$. The morphism $j:G\times \XX\to \IgenX$
  is  an isomorphism.
\end{lem}


\begin{proof}As the stack $\XX$ is separated, we have that the natural
  morphism $I(\XX)\to \XX$ is proper.  As the projection
  $G\times\XX\to \XX$ is a proper morphism, the morphism $j$ is also a
  proper morphism. Its image contains the substack $I(\TT)=\IgenT$
  which is open and dense in $\IgenX$. We deduce that the morphism $j$
  is birational. As the morphism $j$ is \'etale, it is quasi-finite
  (\cf \cite[Expos\'e I.\S 3]{SGA1}).
   The morphism $j$ is proper hence closed and as its image contains the
  open dense torus, $j$ is surjective.  The
  morphism $j$ is a representable, birational, surjective and
  quasi-finite morphism to the smooth \DM stack $\XX$.  The stacky
  Zariski's main theorem \ref{thm:stack-vers-zariski} finishes the
  proof.

\end{proof}



\section{Canonical toric \DM stacks}\label{sec:canon-toric-orbif}
In \S \ref{subsec:canon-stack-assoc} we define  the canonical smooth \DM stack associated to a variety with finite quotient singularities and we show that a canonical smooth \DM stack
satisfies a universal property (Theorem \ref{thm:unicity,can}). This should be well known, but we include it for the reader's convenience.  

In \S\ref{sec:char-canon-toric}, we characterize the canonical
toric \DM stack via its coarse moduli space.

\subsection{Canonical smooth \DM stacks}
\label{subsec:canon-stack-assoc}
In this subsection, we do not assume that smooth \DM stacks are
toric. First, we define canonical smooth \DM stacks and then we prove
their universal property.

We recall a classical result.

\begin{lem}\label{lem:technical,cano}
  Let $S$ be a smooth variety, and $T$ be an affine scheme. Let
  $S'\subset S$ be an open subvariety such that the complement has
  codimension at least $2$ in $S$. Let $f:S'\to T$ be a morphism. Then
  the morphism $f$ extends uniquely to a morphism $S\to T$.
\end{lem}

\begin{proof} The morphism $f$ corresponds to an algebra homomorphism
$K[T]\to \Gamma(S',\OO_{S'})$. Since the complement has codimension
$2$, the restriction map $\Gamma(S,\OO_S)\to\Gamma(S',\OO_{S'})$ is an
isomorphism.
\end{proof}

\begin{defn}\label{defi:codimension,preserving} 
  \begin{enumerate}
  \item 
A dominant morphism $f:V\to W$ of irreducible varieties is called
\textit{codimension preserving} if, for any irreducible closed subvariety
$Z$ of $W$ and every irreducible component $Z_V$ of $f^{-1}(Z)$, one
has $\codim_VZ_V=\codim_WZ$.  
\item A dominant morphism of orbifolds is called
\textit{codimension preserving} if  the induced morphism on every
irreducible component of the coarse moduli spaces is codimension preserving.
\end{enumerate}
\end{defn}

\begin{rem} For any \DM stack, the structure morphism to the coarse
moduli space is codimension preserving. Every flat morphism and in
particular every smooth and \'etale morphism is codimension
preserving. A composition of codimension preserving morphisms is
codimension preserving.
\end{rem}

\begin{defn}\label{defi:canonical}
  Let $\XX$ be an irreducible $d$-dimensional smooth \DM stack. Let
  $\varepsilon:\XX\to X$ be the structure morphism to the coarse
  moduli space. The \DM stack
  $\XX$ will be called \textit{canonical} if the locus where $\eps$ is
  not an isomorphism has dimension $\le d-2$. 
\end{defn}

\begin{rem}\label{rem:cano}
Let $\XX$ be a smooth canonical stack 
  \begin{enumerate}
  \item\label{item:5}  The locus where the
    structure map to the coarse moduli space $\eps:\XX\to X$ is an
    isomorphism is precisely $\eps^{-1}(X_{\sm})$, where $X_{\sm}$ is
    the smooth locus of $X$.
  \item\label{item:1}  The composition of the following isomorphisms 
    \begin{displaymath}
      A^{1}(X)\stackrel{\simeq}{\to} A^{1}(X_{\sm})\stackrel{\simeq}{\to}
      \Pic(X_{\sm})\stackrel{\simeq}{\to}\Pic(\varepsilon^{-1}(X_{\sm}))
\stackrel{\simeq}{\to} \Pic(\XX)
    \end{displaymath} 
is the map sending $[D]$ to $\mathcal{O}(\varepsilon^{-1}(D))$.
    \end{enumerate}
\end{rem}

\begin{thm}[Universal property of canonical smooth \DM stacks]\label{thm:unicity,can} Let $\YY$ be a canonical smooth \DM stack,
  $\eps:\YY\to Y$ the morphism to the coarse moduli space, and
  $f:\XX\to Y$ a dominant codimension preserving morphism with $\XX$ an
  orbifold. Then there exists a unique, up to a unique $2$-arrow, $g:\XX\to\YY$  such that the
  following diagram is commutative
\begin{displaymath}
  \xymatrix{\XX\ar[r]^{\exists\k !\k g}\ar[rd]_{f}& \YY \ar[d]^{\varepsilon} \\ & Y}
\end{displaymath}
\end{thm}

 \begin{proof} 
   We first prove uniqueness. Any two morphisms $g,\bar g$ making the diagram
   commute must agree on the open dense subscheme $f^{-1}(Y_{\sm})$. Put
   $\iota:f^{-1}(Y_{\sm})\hookrightarrow \XX$. Since
   $\YY$ is assumed to be separated, by Proposition \ref{prop:sep,bis}, there
   exists a unique $\alpha: g\to\bar g$ such that $\alpha\ast\Id_{\iota}=\Id$.
   
   By uniqueness, it is enough to prove the result \'etale locally in
   $\YY$, so we can assume that $\YY=[V/G]$ where $V$ is a smooth
   affine variety and $G$ a finite group acting on $V$ without
   pseudo-reflections.
   It is enough to show that there exists an \'etale surjective
   morphism $p:U\to \XX$ with $U$ a smooth variety and a morphism
   $\bar g:U\to \YY$ such that $f\circ p=\eps\circ \bar g$.
 In fact,
   $g$ is defined from $\bar g$ by descent, with the appropriate
   compatibility conditions being taken care of by the uniqueness
   part.
   In this case $Y=V/G$, and $Y_{0}:=V_0/G$ where $V_0\subset V$ is
   the open locus where $G$ acts freely. Let $U_0:=(f\circ
   p)^{-1}(Y_0)$.  As $[V_{0}/G]$ is isomorphic to $Y_{0}$, we have a
   natural morphism $U_{0}\to [V_{0}/G]$. This morphism defines a
     principal $G$-bundle $P_{0}$ on $U_{0}$ and a $G$-equivariant
     morphism $s_{0}:P_{0}\to V_{0}$.
   \begin{equation}\label{eq:12}
     \xymatrix{
P_{0}\ar[rr]^-{s_{0}}\ar[dd]\ar[rd]&&V_{0}\ar'[d][dd]\ar[rd]&\\
&P\ar[dd]\ar[rr]^(.3){s}&&V\ar[dd]\\
U_{0}\ar'[r][rr]\ar[rd]\ar[ddr]&&[V_{0}/G]\ar[rd]\ar[ldd]&\\
&U\ar@{-->}[rr]_-{(P,s)}\ar[ddr]^{f\circ p}&&[V/G]\ar[ddl]^{\varepsilon}\\
&Y_{0}=V_{0}/G\ar[rd]&&\\
&&Y=V/G&\\}
   \end{equation}
   Since
   the $U\setminus U_{0}$ has codimension $\ge 2$, the principal $G$-bundle
   $P_0$ extends uniquely to a principal $G$-bundle $P$ over $U$, and
   by Lemma \ref{lem:technical,cano} (since $V$ is affine) the
   $G$-equivariant morphism $s_{0}:P_0\to V_0$ extends to a morphism
   $s:P\to V$ which is again $G$-equivariant, yielding a morphism
   $\bar g : U\to[V/G]$. The construction above is summarized in the
   $2$-commutative Diagram (\ref{eq:12}) where the squares are
   $2$-cartesian.
 \end{proof}

\begin{cor} Let $\XX$ (resp. $\YY$) be a canonical smooth \DM stack
  with coarse moduli space $X$ (resp. $Y$). Let $\bar f:X\to Y$ be an
  isomorphism. Then there is a unique isomorphism $f:\XX\to \YY$
  inducing $\bar f$.
\end{cor}

 \begin{proof} It is enough to apply the Theorem twice, reversing the role of $\XX$ and $\YY$. 
\end{proof}

\begin{rem}\label{Xcan}
  One can use the corollary to prove the classical fact that every variety $Y$ with
  finite quotient singularities is the coarse moduli space of a
  canonical smooth \DM stack unique up to rigid isomorphism, which we denote by $\YY^{\can}$ (do it
  \'etale locally and then glue).  If $Y$ is the geometric quotient
  $Z/G$ where $Z$ is a smooth variety and $G$ is a group without
  pseudo-reflections acting with finite stabilizers, then
  $\YY^{\can}=[Z/G]$. Notice that this is the case of simplicial toric
  varieties (\cf Section \ref{subsec:geometry-toric}).
\end{rem}

We finish this section with a corollary that will play an important role.
\begin{cor}\label{fmap}
  Let $\XX$ be a  smooth \DM stack with coarse moduli space $\eps:\XX\to
  X$. There is a unique morphism $\XX\to \XX^{\can}$ through which
  $\eps$ factors.
\end{cor}

\begin{proof} 
Apply the theorem with $Y=X$, $\YY=\XX^{\can}$ and $f=\eps$.
\end{proof}

\subsection{The canonical stack of a simplicial toric variety}\label{sec:char-canon-toric}

In this section, we study the canonical stack associated to a
simplicial toric variety.

The main result of this section is the following theorem.



\begin{thm}\label{thm:charac,cano,orb}
  Let $X$ be a simplicial toric variety with torus $T$. Its canonical stack
  $\XX^{\can}$ has a natural structure of toric orbifold such that the
  action  $a^{\can}:T\times \XX^{\can} \to \XX^{\can}$ lifts the
  action $\overline{a}: T\times X \to X$. 
\end{thm}

\begin{proof}
  Denote by $\Sigma$ the fan in $N\otimes_{\zz} \qq$ of the toric variety $X$.
  Without lost of generality, we can assume that the rays of $\Sigma$ generate
  $N\otimes_{\zz} \qq$, so that $X=Z_{\Sigma}/G_A$ (\cf \S
  \ref{subsec:geometry-toric}). The subvariety of points where $G_A$ acts with
  nontrivial stabilizers has codimension $\ge 2$.  Remark \ref{Xcan} implies
  that the canonical stack $\XX^{\can}$ is isomorphic to $[Z_{\Sigma}/G_A]$.  Let
  $T:=(\cc^{\ast})^{n}/G_{A}$ be the torus of the toric variety $X$. Notice that
  $\TT^{\can}=[(\cc^{\ast})^{n}/G_{A}]$ is open dense and isomorphic via
  $\varepsilon\vert_{\TT^{\can}}$ to $T$.  Proposition \ref{prop:action,morphism}
  and the universal property (see Theorem \ref{thm:unicity,can}) of the
  canonical stack imply that the action of $T$ on $X$ lifts to an action of $T$
  on $\XX^{\can}$.
\end{proof}


\begin{rem}\label{rem,can,pic}
  \begin{enumerate}
  \item\label{item:21} Under the hypothesis of Theorem \ref{thm:charac,cano,orb}, we
  have that the restriction of the structure morphism $\varepsilon
  :\XX^{\can}\to X$ to $\TT^{\can}$ is an isomorphism with $T$.

\item\label{item:22}
  Let $\XX$ be a canonical toric \DM stack with \DM torus $\TT=T$ with coarse
  moduli space the simplicial toric variety $X$.  Denote $\Sigma\subset
  N_{\qq}:=N\otimes \qq$ the fan of $X$. Assume that the rays of $\Sigma$
  generates $N_{\qq}$. The proof above shows that $\XX=[Z_{\Sigma}/G_{A}]$ where
  $G_{A}=\Hom(A^{1}(X),\cc^{\ast})=\Hom(\Pic(\XX),\cc^{\ast})$ (\cf Remark
  \ref{rem:cano}.(\ref{item:1})).
\end{enumerate}
\end{rem}

\begin{cor}\label{cor:cano,orbi}Let $\XX$ be a canonical
  toric \DM stack with torus $\mathcal{T}=T$ and coarse
  moduli space the simplicial toric variety $X$.  Denote $\Sigma\subset
  N_{\qq}$ the fan of $X$.
  \begin{enumerate}
  \item\label{item:7} The boundary divisor $\XX\setminus T$ is a simple normal
    crossing divisor, with irreducible components, denoted by $\DD_i$.
 Moreover,
    if  the rays of $\Sigma$ generates $N_{\qq}$, then the divisor $\DD_{i}$
    is isomorphic to $[Z_{i}/G_{A}]$ where $Z_{i}=\{x_{i}=0\}\cap Z_{\Sigma}$.
  \item\label{item:8} The composition morphism  $L\to
    A^{1}(X)\stackrel{\varepsilon^{\ast}}{\to}\Pic(\XX)$ sends $e_{i}$
    to $\OO_{\XX}(\DD_i)$.
  \end{enumerate}
\end{cor}

\begin{proof}
  The first point of the Corollary follows from the fact that the
  inverse image inside $Z_{\Sigma}$ of the torus
  $T=(\cc^{\ast})^{n}/G_{A}$ is $(\cc^{\ast})^{n}$.
  
  The second part of the Corollary follows from the exact sequence
  (\ref{eq:9}) and Remark \ref{rem:cano}.(\ref{item:1}).
\end{proof}

\begin{rem}Let $\XX$ be a canonical toric \DM stack with coarse moduli
    space $X$.
  \begin{enumerate}
  \item Denote by $\Sigma$ the fan of $X$ in $N_{\qq}$. If the
    rays of $\Sigma$ span $N_{\qq}$, from the Corollary and the
    exact sequence (\ref{eq:9}), we have the following exact sequence
    \begin{displaymath}
 \xymatrix{0\ar[r]&M\ar[r]&L\ar[r]&\Pic(\mathcal{X})
\ar[r]&0.}
    \end{displaymath}
  \item For any $i\in\{1, \ldots ,n\}$, the divisor $\mathcal{D}_{i}$
    is Cartier. Hence it corresponds to the invertible sheaf
    $\mathcal{O}(\mathcal{D}_{i})$ with the canonical section $s_{i}$.
    Using Remark \ref{rem:sheaves,global,quotient}, the invertible
    sheaf $\mathcal{O}(\mathcal{D}_{i})$ is associated to the
    representation $G_{A}\to
    G_{L}=(\cc^{\ast})^{n}\stackrel{p_{i}}{\to} \cc^{\ast}$ where
    $p_{i}$ is the $i$-th projection. Moreover, the canonical section
    $s_{i}$ is the $i$-th coordinate of $Z_{\Sigma}$.

  \item Let $\XX$ be a canonical toric \DM stack, then all divisors multiplicities of $\XX$ are equal to $1$ (for the definition of divisor multiplicity see Remark \ref{rem:mult}).
\end{enumerate}
\end{rem}


\section{Toric orbifolds}\label{sec:reduc-toric-orbif}
In this section, we only consider toric \DM stacks with trivial generic
stabilizer that is toric orbifolds.



Let $\XX$ be a smooth \DM stack with coarse moduli space $X$. By
Proposition \ref{prop:toric,orbifold,cms,toric,variety} and Theorem
\ref{thm:charac,cano,orb}, the canonical stack $\XX^{\can}$ has an
induced structure of toric orbifold.  Denote by $\varepsilon_{\XX}:\XX\to
X$ (resp. $\varepsilon_{\XX^{\can}}:\XX^{\can}\to X$) the morphism to the
coarse moduli space.  Theorem \ref{thm:unicity,can} implies that there
exists a unique $f: \XX\to \XX^{\can}$ such that
$\varepsilon_{\XX^{\can}}\circ f= \varepsilon_{\XX}$.

\begin{prop}\label{prop:morph,orb,can,toric} 
  Let $\XX$ be a toric orbifold with torus $T$ and
  coarse moduli space $X$. The canonical morphism $f:\XX\to \XX^{\can}$
  is a morphism of toric \DM stacks where $\XX^{\can}$ is endowed with
  the induced structure of toric orbifold.
\end{prop}

 \begin{proof}
   The universal property of
   the canonical stack (\cf Theorem \ref{thm:unicity,can}) applied to
   $\Id:T \to T$ implies that $f\vert_{T}:T\to \TT^{\can}$.  

 \end{proof}
 Notice that the morphism $f\vert_{T}:T \to \TT^{\can}$ in the proof
 above is an isomorphism because $\XX$ is a toric  orbifold. 
 
 Denote
$\bs{D}^{\can}:=(D_{1}^{\can}, \ldots ,D_{n}^{\can})$ (\cf \ref{subsubsec:root-line-section}).

\begin{thm}\label{thm:charac,reduced,toric}
  \begin{enumerate}
  \item\label{item:2} Let $X$ be a simplicial toric variety with torus
    $T$. Denote by $\Sigma$ a fan of $X$. For each ray $\rho_{i}$ of $\Sigma$,
    choose $a_{i}$ in $\nn_{>0}$. Denote $\bs{a}:=(a_{1}, \ldots
    ,a_{n})\in (\nn_{>0})^{n}$.
    Then $\sqrt[\bs{a}]{\bs{D}^{\can}/\XX^{\can}}$ has a unique
    structure of toric orbifold with torus $T$ such that the canonical
    morphism $\pi:\sqrt[\bs{a}]{\bs{D}^{\can}/\XX^{\can}}\to \XX^{\can}$ is a
    morphism of toric \DM stacks with divisor multiplicities $\bs{a}$.
  \item\label{item:6} Let $\XX$ be a toric orbifold with coarse moduli
    space $X$. Let $\bs{a}:=(a_{1}, \ldots ,a_{n})$ be its divisors
    multiplicities.  Then $\XX$ is naturally isomorphic as toric \DM
    stack to $\sqrt[\bs{a}]{\bs{D}^{\can}/\XX^{\can}}$ defined in (1).
\end{enumerate}
\end{thm}



\begin{proof}
  (1) Let $\TT^{\can}\subseteq \XX^{\can}$ be the inverse image of $T$ (which is
  isomorphic to $T$). Note that
  $\pi^{-1}(T)\subseteq\sqrt[\bs{a}]{\bs{D}^{\can}/\XX^{\can}}$ is isomorphic to
  $\TT^{\can}$ by Property 2 of Section \ref{subsubsec:root-line-section}. Let $j:
  T\to\sqrt[\bs{a}]{\bs{D}^{\can}/\XX^{\can}}$ be the dominant open embedding.
  We need to prove that $T$ acts on $\sqrt[\bs{a}]{\bs{D}^{\can}/\XX^{\can}}$
  compatibly with $j$. We know that $T$ acts on $\XX^{\can}$. To define
  $T\times\sqrt[\bs{a}]{\bs{D}^{\can}/\XX^{\can}}\to
  \sqrt[\bs{a}]{\bs{D}^{\can}/\XX^{\can}} $ we use the universal property and
  the fact that $\bs{D}^{\can}\subseteq \XX^{\can}$ is $T$-invariant.
  
    (2) For any
    $i\in\{1, \ldots ,n\}$, denote by $D_{i},D_{i}^{\can},\DD_{i}(\XX)$
    the divisor corresponding to the ray $\rho_{i}$ in respectively $X,
    \XX^{\can}$ and $\XX$.  Theorem \ref{thm:charac,cano,orb}
    implies there exist a unique morphism $f:\mathcal{X}\to \XX^{\can}$
    such that $\varepsilon_{\XX^{\can}}\circ f= \varepsilon_{\XX}$.  By
    definition of the divisors multiplicities, for any ray $\rho_{i}$,
    we have $f^{-1}D^{\can}_{i}=a_{i}\mathcal{D}_{i}(\XX)$. The
    Cartier divisors $\bs{\DD}(\XX):=(\mathcal{D}_{1}(\XX), \ldots
    ,\DD_{n}(\XX))$ define a morphism $\mathcal{X}\to
    [\aaa^{n}/(\cc^{\ast})^{n}]$ such that the following diagram is
    $2$-commutative :
    \begin{equation}\label{eq:16}
\xymatrix{\XX \ar[rr]^-{\bs{\DD}(\XX)}\ar[d]^{f}&& 
\left[\aaa^{n}/(\ccs)^{n}\right] \ar[d]^{\wedge\bs{a}}\\ 
\XX^{\can}\ar[rr]^-{\bs{D}^{\can}}&&\left[\aaa^{n}/(\ccs)^{n}\right]}
    \end{equation}
where the morphism $\wedge \bs{a}$ is defined in Section
\ref{subsubsec:root-line-section}.
By the universal property of fiber product, we deduce a unique
morphism $g:\mathcal{X}\to \sqrt[\bs{a}]{\bs{D}^{\can}/\XX^{\can}}$
such that the following diagram is strictly commutative
   \begin{displaymath}
     \xymatrix{\mathcal{X}\ar[r]^-{g}\ar[rd]_{f}
     &\sqrt[\bs{a}]{\bs{D}^{\can}/\XX^{\can}}\ar[d]^{\pi}\\ &\XX^{\can}}
   \end{displaymath}


   We will use the Zariski's main theorem (\cf Theorem
   \ref{thm:stack-vers-zariski} ) to prove that $g$ is an isomorphism.
   We first notice that $\sqrt[\bs{a}]{\bs{D}^{\can}/\XX^{\can}}$ is
   smooth for property (\ref{item:9}) in
   \S\ref{subsubsec:root-line-section}. As, the restriction of $g$
   over $\mathcal{X}^{\can}-\cup_{i,j} D_{i}^{\can}\cap D_{j}^{\can}$
   is an isomorphism, the morphism $g$ is birational. Notice that
   $\cup_{i,j} D_{i}^{\can}\cap D_{j}^{\can}$ is a subset of
   codimension $\geq$ 2. The morphism $g$ is proper, hence closed, so
   we deduce that $g$ is also surjective because its image contains
   the dense torus. Let us show that $g$ is representable and \'etale.
   Let $S$ be a scheme. Consider the following $2$-cartesian diagram
   \begin{displaymath}
     \xymatrix@1{\mathcal{Y}\ar@{}[rrd]|{\square}\ar[rr]^-{\overline{g}}\ar[d]&&S\ar[d]\\\mathcal{X}\ar[rr]^-{g}&&\sqrt[\bs{a}]{\bs{D}^{\can}/\XX^{\can}}}
   \end{displaymath}
   Let $U\to \mathcal{Y}$ be an \'etale atlas of $\mathcal{Y}$. First
   we observe that the morphism $U\to S$, denote it by
   $\widetilde{g}$, must be flat, so that the morphism $g$ is flat
   too. To verify this we can apply \cite[Thm 23.1]{Mca} using that
   both $S$ and $U$ are smooth and the dimension of the fibers of
   $\widetilde{g}$ is constantly zero. To prove that the dimension of
   the fibers is zero we just need to observe that both $\pi$ and $f$
   are quasi-finite, since they are morphisms from a stack to its
   coarse moduli space, and $f$ factors through $g$ so that it must be
   quasi-finite too.  We now note that the morphism $U\to S$ is
   \'etale away from a codimension $\geq$ 2 subset, so we can apply
   the theorem of purity of branch locus (\cf Theorem 6.8 p.125 in
   \cite{Altman-Kleiman-Intro-Grothendieck-duality}) and deduce that
   $U\to S$ is \'etale \ie $\overline{g}:\mathcal{Y}\to S$ is \'etale.
   Without loss of generality we can assume that $S$ is actually an
   atlas; we assume that $\mathcal{Y}$ is a stack and we prove that it
   must be actually a scheme. First of all we observe that it cannot
   have generically non trivial stabilizer, since the morphism $\YY\to
   \XX$ is representable it must induce an injection of the stabilizer
   at each geometric point
   \cite{Abramovich-Vistoli-Compactification-stable-maps}, but $\XX$
   is an orbifold so that $\YY$ must be an orbifold too. There exists
   an \'etale representable map $[V/K]\to \mathcal{Y}$ where $V$ is a
   smooth variety and $K$ is a finite group. Hence the induced map
   $V\to S$ is \'etale. By the universal property, it factors via the
   coarse moduli space $V/K$, and the map $V\to V/K$ is not injective
   on tangent vectors unless $K$ is acting freely, hence $V\to V/K$
   cannot be \'etale unless $\mathcal{Y}$ has trivial stabilizers
   everywhere. We now observe that the morphism $V/K\to S$ is still
   birational surjective and quasi-finite, using Zariski's main
   theorem for schemes we can deduce that it is an isomorphism, in
   particular it is \'etale and this implies that $V\to V/K$ must be
   \'etale.  We conclude that $\mathcal{Y}$ is a scheme \ie $g$ is
   representable and \'etale.  So it is also quasi-finite (\cf
   \cite[Expos\'e I.\S 3]{SGA1}).

   As the morphism $g$ is representable, surjective, birational and
   quasi-finite, the stacky Zariski's main theorem
   \ref{thm:stack-vers-zariski} implies that $g$ is an isomorphism.

\end{proof}

The following Corollary is a consequence of Property \ref{item:9}   of
Section \ref{subsubsec:root-line-section} and Theorem \ref{thm:charac,reduced,toric}.
\begin{cor}\label{cor:reduced,toric}
  Let $\mathcal{X}$ be a toric orbifold with coarse moduli space $X$.
  The reduced closed substack $\mathcal{X}\setminus \mathcal{T}$ is a simple
  normal crossing divisor.
\end{cor}

\begin{rem}\label{rem:picard,group,orb} Let $\XX$ be a toric
    orbifold with coarse moduli space $X$.  Diagram (\ref{eq:29}) and
    Theorem \ref{thm:charac,reduced,toric} imply that we have the
    following morphism of exact sequences
\begin{equation}\label{eq:28}
\xymatrix{0 \ar[r]&\zz^{n}\ar[d]\ar[r]^{\times \bs{a}}&\zz^{n}\ar[r]\ar[d]&
  \oplus_{i=1}^{n} \zz/a_{i}\zz \ar[r]\ar@{=}[d]&0
\\0\ar[r]& \mathbb\Pic(\XX^{\can})\ar[r]^-{f^{\ast}} &\Pic(\mathcal{X})\ar[r]&\oplus_{i=1}^{n} \zz/a_{i}\zz\ar[r]&0}
\end{equation}
where the vertical morphisms sends $1\mapsto
  \OO(D_{i}^{\can})$ and $1\mapsto \OO(\DD_{i})$. 
\end{rem}


\renewcommand\GG{\Gm}

\section{Toric \DM stacks}\label{sec:toric-dm-stacks-1}
In this section we will show that each toric \DM stack is isomorphic to a fibered product of
root stacks on its rigidification.  To prove this theorem, we will
recall in Section \ref{subsec:gerbe+root} the relation between banded
gerbes and root constructions. Then we will show in Theorem
\ref{thm:gerbe,banale} that any toric \DM stack is an essentially
trivial gerbe on its rigidification.  In Section
\ref{sec:char-toric-dm}, we will prove the main result in Theorem
\ref{thm:cara,gerbe,toric}.
\subsection{Gerbes and root constructions}
\label{subsec:gerbe+root}

First, we recall some general notion on banded gerbes (\textit{gerbes
  li\'ees}). We refer to \cite{Gcna-1971} chapter IV.2 for a complete
treatment and to Section 3 of \cite{EKVbgqs-2001} for a shorter
reference. Let $\XX$ be a smooth \DM stack. Let $G$ be an abelian sheaf
of groups\footnotemark \footnotetext{The non abelian case has a richer
  structure but for the sake of simplicity we just skip all these
  additional features and refer the interested reader to
  \cite{Gcna-1971}.} and $\mathcal{G}\to\XX$ a gerbe.  For every
\'etale chart $U$ of $\XX$ and every object $x\in\mathcal{G}(U)$ let
$\alpha_x:G\vert_U\to\Aut_U(x)$ be an isomorphism of sheaves of groups
such that the natural compatibilities coming from the fibered
structure of the gerbe are satisfied. The collection of these
isomorphisms is called a $G$-\textit{banding}. A $G$-\textit{banded} gerbe is the data
of a gerbe and a $G$-banding. Two $G$-banded gerbes are said to be
$G$-\textit{equivalent} if they are isomorphic as stacks and the isomorphism
makes the two bandings compatible in the natural way. Giraud proved in
\cite{Gcna-1971} (chapter IV.3.4) that the group $H_{\et}^2(\XX,G)$
classifies equivalence classes of $G$-banded gerbes.

\begin{rem}\label{rem:observations}
  We anticipate some observations about the banding which will be
  useful in the following:
\begin{enumerate}
\item\label{item:15} The $b$-th root of a line bundle on $\XX$ is a gerbe which is
  banded in a natural way by the constant sheaf $\mu_b$; the banding
  is the canonical isomorphism between the group of automorphisms of
  any object and $\mu_b$.
\item\label{item:16} Given $\mathcal{G}\to \XX$ a $G$-banded gerbe, every
  rigidification of $\mathcal{G}$ by a subgroup $H$ of $G$ inherits a
  $(G/H)$-banding from the $G$-banding of $\mathcal{G}$.
\end{enumerate}   
\end{rem}

Here we introduce the concept of an essentially trivial gerbe which will play an important role in this section.
The Kummer sequence 
\begin{displaymath}
  \xymatrix{1\ar[r]&\mu_{b}\ar[r]^{\iota}&\Gm\ar[r]^{\wedge b}&\Gm\ar[r]&1}
\end{displaymath}
induces  the long exact sequence 
\begin{equation}\label{eq:17}
  \xymatrix{\cdots\ar[r]&H^{1}_{\et}(\XX,\Gm)\ar[r]^{\partial}&H^{2}_{\et}(\XX,\mu_{b})\ar[r]^{\iota_{\ast}}&H^{2}(\XX,\Gm)\ar[r]&\cdots}
\end{equation}

\begin{defn}
  A $\mu_{b}$-banded gerbe in $H^{2}_{\et}(\XX,\mu_{b})$ is \textit{essentially trivial} if its image by $\iota_{\ast}$ is the trivial gerbe in
  $H^{2}_{\et}(\XX,\Gm)$.
\end{defn}

\begin{rem}\label{rem:gerbe}
  \begin{enumerate}
  \item\label{item:10} It follows from Section \ref{sec:root-line} that a
    $\mu_{b}$-banded gerbe is essentially trivial if and only if it
    is a $b$-th root of an invertible sheaf on $\XX$.
  \item\label{item:14} As the $\mu_{b}$-banded gerbe
    $\sqrt[b]{L\otimes M^{\otimes b}/\XX}$ is isomorphic to
    $\sqrt[b]{L/\XX}$, we deduce a bijection between essentially
    trivial $\mu_{b}$-banded gerbes and $\Pic(\XX)/b\Pic(\XX)$.
\end{enumerate}
\end{rem}

\begin{lem}\label{lem:ess,trivial,gerb}
 There is a natural bijection between essentially trivial gerbes in
  $H^{2}_{\et}(\XX,\mu_{b})$ and elements in $\Ext^{1}(\zz/b\zz,\Pic(\XX))$.
\end{lem}

\begin{proof}By Remark \ref{rem:gerbe}.(\ref{item:14}), it is enough
  to show that $ \Ext^{1}(\zz/b\zz,\Pic(\XX))$ is isomorphic to
  $\Pic(\XX)/b\Pic(\XX)$. This follows from  the exact sequence
  \begin{displaymath}
    \xymatrix{\Hom(\zz,\Pic(\XX))\ar[r]^-{\wedge b}&\Hom(\zz,\Pic(\XX))\ar[r]&\Ext^{1}(\zz/b,\Pic(\XX))\ar[r]&0}
  \end{displaymath}

\end{proof}

  Let $G$ be a finite abelian group. Fix a decomposition 
  $G=\prod_{j=1}^{\ell}\mu_{b_{j}}$. We deduce an isomorphism 
  \begin{align}\label{eq:19}
    H^{2}_{\et}(\XX,G)&\longrightarrow
    \bigoplus_{j=1}^{\ell}H^{2}_{\et}(\XX,\mu_{b_{j}}) \\
\alpha & \longmapsto (\alpha_{1}, \ldots ,\alpha_{\ell}) \notag
  \end{align}

\begin{defn}
  Let $G$ be a finite abelian group. A $G$-banded gerbe associated to
  $\alpha\in H^{2}(\XX,G)$ is \textit{essentially trivial} if there is a
  decomposition of $G=\prod_{j=1}^{\ell}\mu_{b_{j}}$ such that for any
  $j\in\{1, \ldots ,\ell\}$, the $\mu_{b_{j}}$-banded gerbe $\alpha_{j}$
  is essentially trivial.
\end{defn}

\begin{rem} Being essentially trivial does not depend on the choice of a decomposition of $G$.
\end{rem}

\begin{prop}\label{prop:ess,trivi,gerb} Let $G$ be a finite abelian group. Fix a decomposition of
  $G=\prod_{j=1}^{\ell}\mu_{b_{j}}$.
  There are bijections between
  \begin{align*}
    &\left\{\mbox{Essentially trivial gerbes in }
    \oplus_{j=1}^{\ell}H^{2}_{\et}(\XX,\mu_{b_{j}})\right\} \\
\stackrel{1:1}{\longleftrightarrow}& \left\{\mbox{Fibered products over }
    \XX \mbox{ of }b_{j}\mbox{-th roots of invertible
    sheaves}\right\}\\ \stackrel{1:1}{\longleftrightarrow}&
    \prod_{j=1}^{\ell}\Pic(\XX)/b_{j}\Pic(\XX)
\stackrel{1:1}{\longleftrightarrow}
    \prod_{j=1}^{\ell}\Ext^{1}(\zz/b_{j}\zz,\Pic(\XX))
  \end{align*}
\end{prop}

\begin{rem}
  To be more concrete, let us explicitly describe the last bijection. For the
  sake of simplicity, we consider the case $j=1$.
  To the class $[L_{0}]$ in 
  $\Pic(\XX)/b\Pic(\XX)$, we associate the
  extension
  \begin{displaymath}
    \xymatrix{0\ar[r]&\Pic(\XX)\ar[r]&\Pic(\XX)\times_{\Pic(\XX)/b\Pic(\XX)}\zz/b\zz\ar[r]& \zz/b\zz \ar[r]&0}
  \end{displaymath} where the fiber product is given by the standard
  projection $\Pic(\XX)\to\Pic(\XX)/b\Pic(\XX)$ and the morphism
  $\zz/b\zz\to \Pic(\XX)$ that sends the class of $1$ to the class
  $[L_{0}]$. The first morphism in the extension sends the invertible
  sheaf $L$ to $(L^{\otimes b},0)$.
  
  Let $0\to\Pic(\XX)\to A \to \zz/b\to 0$ be an extension. We consider
  the projective resolution $0\to\zz\stackrel{\times b}{\to}\zz\to
  \zz/b\to 0$.  There exists $f$ and $\widetilde{f}$ such that the
  following diagram is a morphism of short exact sequences.
   \begin{displaymath}
     \xymatrix{0\ar[r]&\zz\ar[r]\ar[d]^{\widetilde{f}}&\zz \ar[r]\ar[d]^{f}&\zz/b \ar[r]\ar@{=}[d]& 0 \\ 0\ar[r]&\Pic(\XX)\ar[r]& A\ar[r]& \zz/b\ar[r]&0}
   \end{displaymath}
   The class $[\widetilde{f}(1)]$ in $\Pic(\XX)/b\Pic(\XX)$ is the
   element that corresponds to the above extension. Notice that
   different liftings $f,\widetilde{f}$ lead to different elements in
   $\Pic(\XX)$ with the same class in $\Pic(\XX)/b\Pic(\XX)$.

The two maps define above are inverse to each other.
\end{rem}

\begin{proof}[Proof of Proposition \ref{prop:ess,trivi,gerb}]
  Most of the Proposition is a direct consequence of
  Remark~\ref{rem:gerbe} and Lemma \ref{lem:ess,trivial,gerb}. The
  only non trivial fact  to prove is that an essentially trivial gerbe
  defined by $\alpha=(\alpha_{1}, \ldots ,\alpha_{\ell})\in
  \oplus_{j=1}^{\ell}H^{2}_{\et}(\XX,\mu_{b_{j}})$ is given by a fiber
  product of the gerbes defined by the $\alpha_{j}$'s.  Without loss
  of generality, we can assume that $\alpha=(\alpha_{1},\alpha_{2})$; the
  general case is proved by induction. The gerbe defined by
  $\alpha_{1}$ (resp. $\alpha_2$) is isomorphic to the rigidification
  $\mathcal{G}\fs \mu_{b_{2}}$ (resp. $\mathcal{G}\fs \mu_{b_{1}}$ ).
  Hence we have the following $2$-commutative diagram
  \begin{displaymath}
    \xymatrix{
 & \mathcal{G}\fs \mu_{b_{2}}\ar[dr] & \\
 \mathcal{G}\ar[ur]\ar[rr]\ar[dr] &  & \XX \\
 &\mathcal{G}\fs \mu_{b_{1}}\ar[ur] & \\ 
}
  \end{displaymath}
Remark \ref{rem:observations}.(\ref{item:16}) implies that $\mathcal{G}\to\mathcal{G}\fs \mu_{b_{1}} $ (resp. $\mathcal{G}\to\mathcal{G}\fs \mu_{b_{2}} $ ) is a $\mu_{b_2}$-banded gerbe (resp. $\mu_{b_1}$-banded).
By the universal property of the fiber product we are given a morphism 
$\mathcal{G}\to \mathcal{G}\fs \mu_{b_{1}}\times_{\XX} \mathcal{G}\fs
\mu_{b_{1}}$. Two gerbes banded by the same group over the same base $\XX$ are
either isomorphic as stacks or they have no morphisms at all; this completes the proof.
\end{proof}

\subsection{Gerbes on toric orbifolds}
\label{subsec:gerbe,toric,orb}

\begin{thm}\label{thm:gerbe,banale}
  Let $\mathcal{X}$ be a toric orbifold with torus $T$. Denote
by  $\iota:T\hookrightarrow \mathcal{X}$ the immersion of the 
  torus. Then the morphism
\begin{displaymath}
\iota^{\ast}:H^{2}_{\et}(\mathcal{X},\GG) \to H^{2}_{\et}(T,\GG)
\end{displaymath}
is injective.
\end{thm}

Notice that in the following proof we will use that a toric orbifold
is a global quotient $[Z_{\Sigma}/G_{\XX}]$ where
$G_{\XX}:=\Hom_{\zz}(\Pic(\mathcal{X}),\cc^{\ast})$. This will be
proved in Theorem \ref{thm:toric,global,quotient} and does not depend on the
results of this subsection.

We first proof some preliminary results.
\begin{lem}[Artin]\label{lem:remove-cod-2}
  Let $S$ be a smooth quasi-projective variety, $S_2\subseteq S$
  a closed subscheme of codimension $\geq 2$. Then the natural map
  $H_{\et}^i(S,\GG)\to H_{\et}^i(S\setminus S_2,\GG)$ is an isomorphism for all $i$.
\end{lem}

\begin{proof}
  The statement is obvious if we replace sheaf cohomology with \v{C}ech
  cohomology. To prove the lemma, we just apply Corollary $4.2$ p.295 of 
\cite{Ajhr71} (see also Theorem 2.17 p.104 of \cite{Met}).
\end{proof}

\begin{lem}[Olsson]\label{lem:olsson-spectral}
  Let $\mathcal{X}$ be an Artin stack and $X_0$ an atlas. Denote by
  $X_{p}=X_{0}\times_{\mathcal{X}}\cdots\times_{\mathcal{X}}X_{0}$. Let
  $\mathcal{F}$ be an abelian sheaf of groups on $\mathcal{X}$ and
  $\mathcal{F}_p$ its restriction to $X_{p}$. There is a spectral sequence with
  $E_{1}^{pq}(\mathcal{X}):=H_{\et}^{q}(X_{p},\mathcal{F}_{p})$ that abuts to
  $H_{\et}^{p+q}(\mathcal{X},\mathcal{F})$.
\end{lem}

\begin{proof}
  This lemma follows immediately from Corollary 2.7 p.4 and Theorem 4.7 p.13 in
  \cite{Osas05}.
\end{proof}

 \begin{proof}[Proof of Theorem \ref{thm:gerbe,banale}]
Let $\XX$ be a toric orbifold with coarse moduli space a simplicial toric variety
 $X$. Denote by $\Sigma\subset N_{\qq}$ the fan of $X$. Without lost of
 generality, we can assume that the rays of $\Sigma$ generate $N_{\qq}$.
  
By Theorem \ref{thm:charac,reduced,toric} in the case of orbifolds  and
   Lemma \ref{lem:root}, we have that
   $\mathcal{X}=[Z_{\Sigma}/G_{\XX}]$ where
   $G_{\XX}:=\Hom_{\zz}(\Pic(\mathcal{X}),\cc^{\ast})$. 
Denote by $n$ the number of
   rays of the fan $\Sigma$.
   Put
 \begin{displaymath}
   Z_2:=\{z\in Z_{\Sigma}\subset\cc^{n}\vert \forall i\in\{1, \ldots ,n\},
   \prod_{j\neq i} z_{j}=0\}
 \end{displaymath}
the union  of $T$-orbits in $Z_{\Sigma}$ of codimension $\geq 2$.
The closed subscheme $Z_2$ of $Z_{\Sigma}$ is of codimension $2$. Hence
 the quotient stack $[(Z_{\Sigma}\setminus Z_{2})/G_{\XX}]$ is a closed substack of
 codimension $2$ of $\mathcal{X}$.  For any $i\in\{1, \ldots ,n\}$, put
 \begin{displaymath}
   U_{i}:=\{z\in Z_{\Sigma}\subset\cc^{n}\vert \forall j \in 
\{1, \ldots ,n\}\setminus \{i\}, z_{j}\neq 0\}.
 \end{displaymath}
 We have that $U_{i}$ is isomorphic to $\aaa^{1}\times
 (\cc^{\ast})^{n-1}$ and
 that  the natural morphism
 \begin{align*}
   \coprod_{i\in\{1, \ldots ,n\}} U_{i} \to Z_{\Sigma}\setminus Z_2
 \end{align*}
 is \'etale and surjective.
 We deduce that $\coprod_{i\in\{1, \ldots ,n\}} [U_{i}/G_{\XX}] \to
 [(Z_{\Sigma}\setminus Z)/G_{\XX}]$ is  \'etale and surjective.
 Put $X_{0}:=\coprod_{i\in\{1, \ldots ,n\}}U_{i}$. The natural
 morphism $X_{0} \to [(Z_{\Sigma}\setminus Z_2)/G_{\XX}]$ is an \'etale atlas.
 
 Denote by $X_{p}=X_{0}\times_{\mathcal{X}}\cdots\times_{\mathcal{X}}X_{0}$.
 From Lemma \ref{lem:olsson-spectral} we have a spectral sequence
 $E_{1}^{pq}(\mathcal{X}):=H_{\et}^{q}(X_{p},\GG\vert_{X_{p}})$ abutting to
 $H_{\et}^{p+q}([(Z_{\Sigma}\setminus Z_2)/G_{\XX}],\GG)$. Using this spectral
 sequence and Lemma \ref{lem:remove-cod-2} we obtain that the natural
 morphism $ H_{\et}^{i}(\mathcal{X},\GG)=H_{\et}^{i}([(Z_{\Sigma} \setminus
 Z_2)/G_{\XX}],\GG)$ is an isomorphism for $i=(0,1,2)$.  Finally, the theorem
 follows from Lemmas \ref{lem:exact,seq} and \ref{lem:exact,seq,injective}.
 \end{proof}

\begin{lem}\label{lem:exact,seq}
  We have the following morphism of short exact sequences
  \begin{displaymath}
    \xymatrix{0\ar[r]&E_{4}^{20}(\mathcal{X})\ar[r]\ar[d]^{\alpha}&H_{\et}^{2}(\mathcal{X},\GG)\ar[r]\ar[d]^{j^{\ast}}&E_{2}^{02}(\mathcal{X})\ar[r]\ar[d]^{\beta}&0
    \\0\ar[r]&E_{4}^{20}(\mathcal{T})\ar[r]&H_{\et}^{2}(\mathcal{T},\GG)\ar[r]&E_{2}^{02}(\mathcal{T})\ar[r]&0 }
  \end{displaymath}
\end{lem}

\begin{lem}\label{lem:exact,seq,injective}
The vertical maps $\alpha:E_{4}^{20}(\mathcal{X})\to
E_{4}^{20}(\mathcal{T})$ and $\beta:E_{2}^{02}(\mathcal{X})\to
E_{2}^{02}(\mathcal{T})$ are injective.
\end{lem}

 \begin{proof}[Proof of Lemma \ref{lem:exact,seq}]
 To prove the lemma, we are just interested in
  $E_{\infty}^{pq}(\mathcal{X})$ for $ p+q=2$.
 We start by proving that we have
 \begin{equation}\label{eq:1}
    \xymatrix{0\ar[r]&E_{\infty}^{20}(\mathcal{X})\ar[r]&H_{\et}^{2}(\mathcal{X},\GG)\ar[r]&E_{\infty}^{02}(\mathcal{X})\ar[r]&0}
  \end{equation}
   Hilbert's Theorem 90
  (\cf Proposition 4.9 of \cite{Met}) implies that
     \begin{displaymath}
       H^{1}_{\et}(X_{p},\GG)=H^{1}_{\text{Zariski}}(X_{p},\mathcal{O}_{X_{p}}^{\ast})=\Pic(X_{p}).
     \end{displaymath}
     Using the notation of the proof of Theorem \ref{thm:gerbe,banale}, for any
     ray $i\in \{1, \ldots ,n\}$, we have that $[U_{i}/G_{\mathcal{X}}]$ is
     isomorphic to $[\aaa^{1}/\mu_{a_{i}}]\times (\cc^{\ast})^{n-1}$ where
     $a_{i}$ is the multiplicity along the divisor $\mathcal{D}_{i}$ (\cf
     Remark \ref{rem:mult}). Hence, we have that $ X_{p}=\coprod_{i_{0},
       \ldots ,i_{p}\in\{1, \ldots ,n\}} U_{i_{0}\cdots i_{p}}$ where
 \begin{displaymath}
   U_{i_{0}\cdots i_{p}}=
   \begin{cases}U_{i_{0}}\times \mu_{a_{0}}^{p+1} &
     \mbox{if } i_{0}=\cdots=i_{p}\\
     \mathcal{T}&\mbox{otherwise.}
   \end{cases}
   \end{displaymath}
   Hence, for any $p$ we have that
   $E_{1}^{p1}(\mathcal{X})=E_{\infty}^{p1}(\mathcal{X})=H^{1}_{\et}(X_{p},\GG)=0$.
   We deduce the exact sequence (\ref{eq:1}).

 We now show that
  $E_{\infty}^{20}(\mathcal{X})=E_{4}^{20}(\mathcal{X})$ and
 $E_{\infty}^{02}(\mathcal{X})=E_{2}^{02}(\mathcal{X})$.

  \begin{figure}[ht]
  \begin{center}
  \psfrag{p}{\tiny{$p$}}
  \psfrag{q}{\tiny{$q$}}
  \psfrag{q}{\tiny{$q$}}
  \psfrag{0}{\tiny{$0$}}
  \psfrag{d120}{\tiny{$d_{1}$}}
  \psfrag{d121}{\tiny{$d_{1}$}}
  \psfrag{d122}{\tiny{$d_{1}$}}
  \psfrag{d100}{\tiny{$d_{1}$}}
  \psfrag{d101}{\tiny{$d_{1}$}}
  \psfrag{d102}{\tiny{$d_{1}$}}
  \psfrag{d2}{\tiny{$d_{2}$}}
  \psfrag{d3}{\tiny{$d_{3}$}}
    \includegraphics[width=0.3\linewidth]{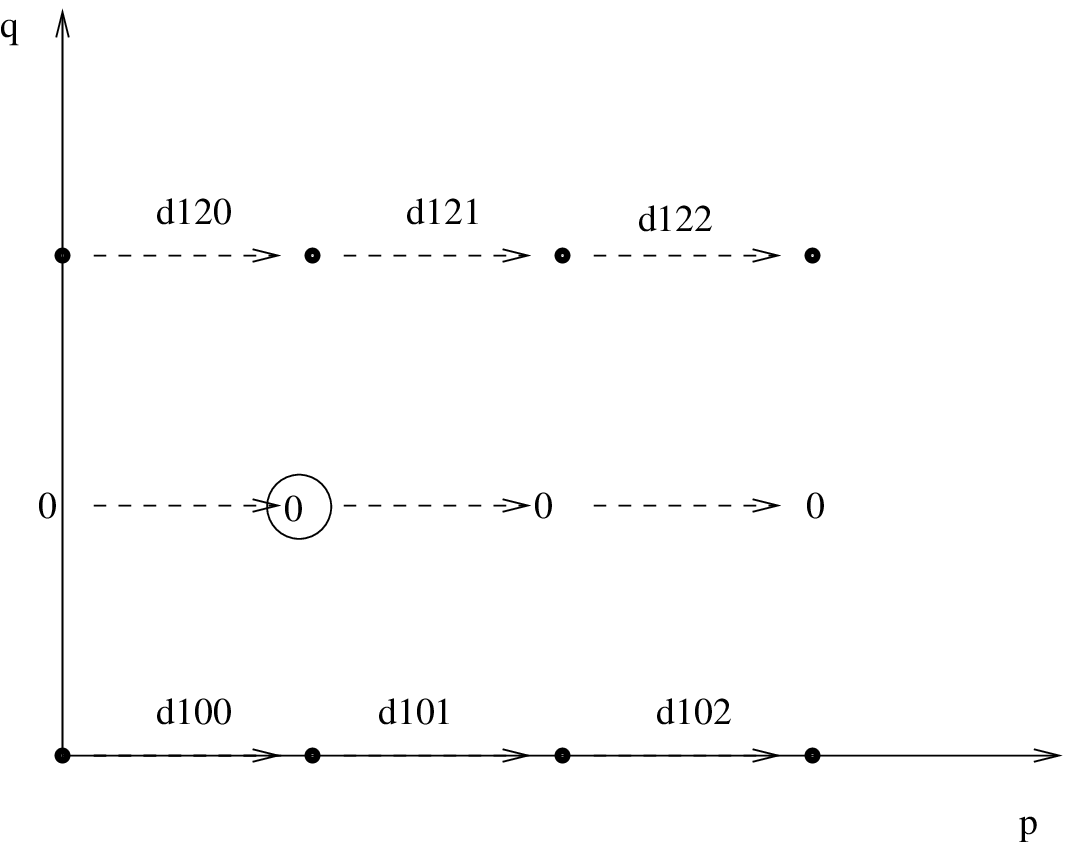}
    \includegraphics[width=0.3\linewidth]{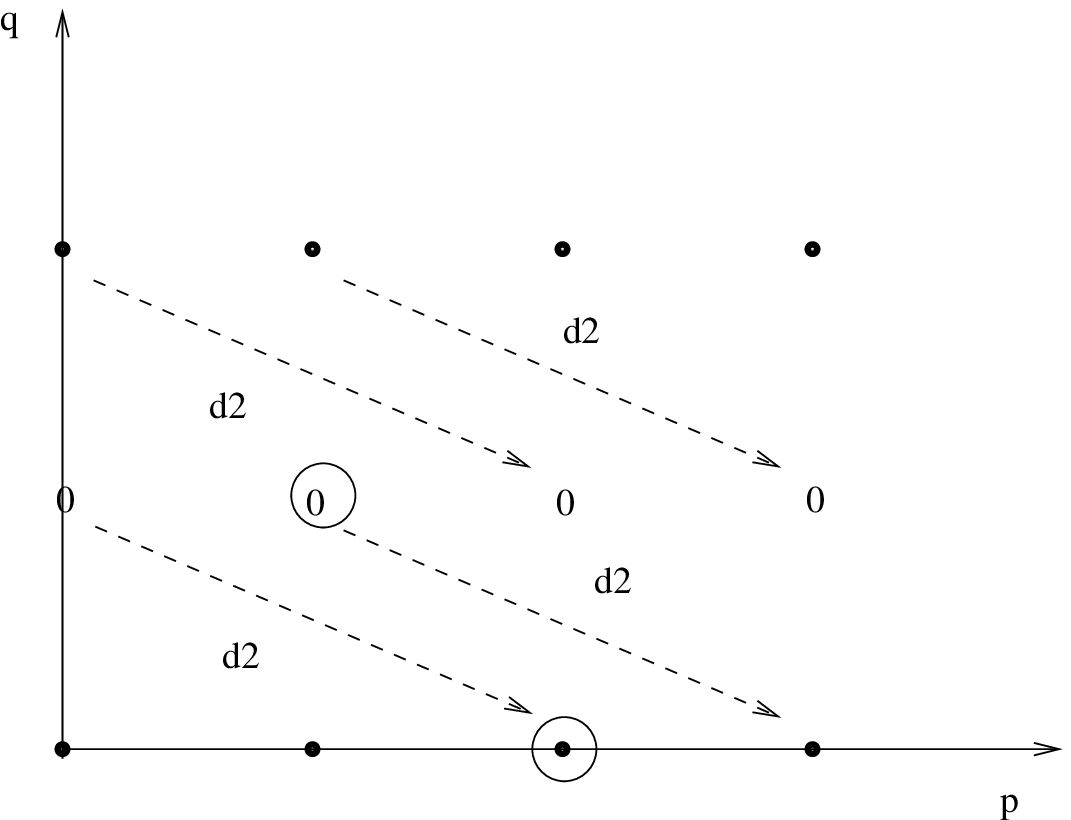}
  \caption{Terms $E^{pq}_{1}(\mathcal{X})$ and $E^{pq}_{2}(\mathcal{X})$}\label{fig:E1,E2}
  \end{center}
  \end{figure}

  \begin{figure}[ht]
  \begin{center}
  \psfrag{p}{\tiny{$p$}}
  \psfrag{q}{\tiny{$q$}}
  \psfrag{q}{\tiny{$q$}}
  \psfrag{0}{\tiny{$0$}}
  \psfrag{d120}{\tiny{$d_{1}$}}
  \psfrag{d121}{\tiny{$d_{1}$}}
  \psfrag{d122}{\tiny{$d_{1}$}}
  \psfrag{d100}{\tiny{$d_{1}$}}
  \psfrag{d101}{\tiny{$d_{1}$}}
  \psfrag{d102}{\tiny{$d_{1}$}}
  \psfrag{d2}{\tiny{$d_{2}$}}
  \psfrag{d3}{\tiny{$d_{3}$}}
    \includegraphics[width=0.3\linewidth]{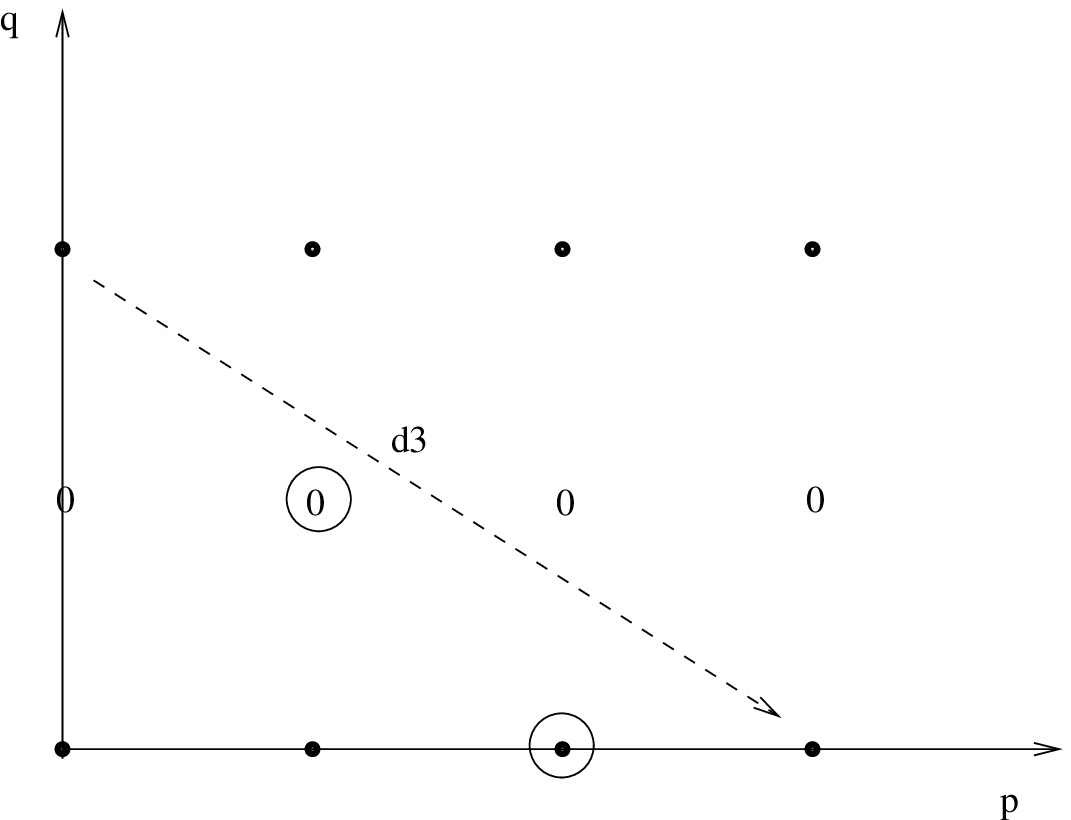}
    \includegraphics[width=0.3\linewidth]{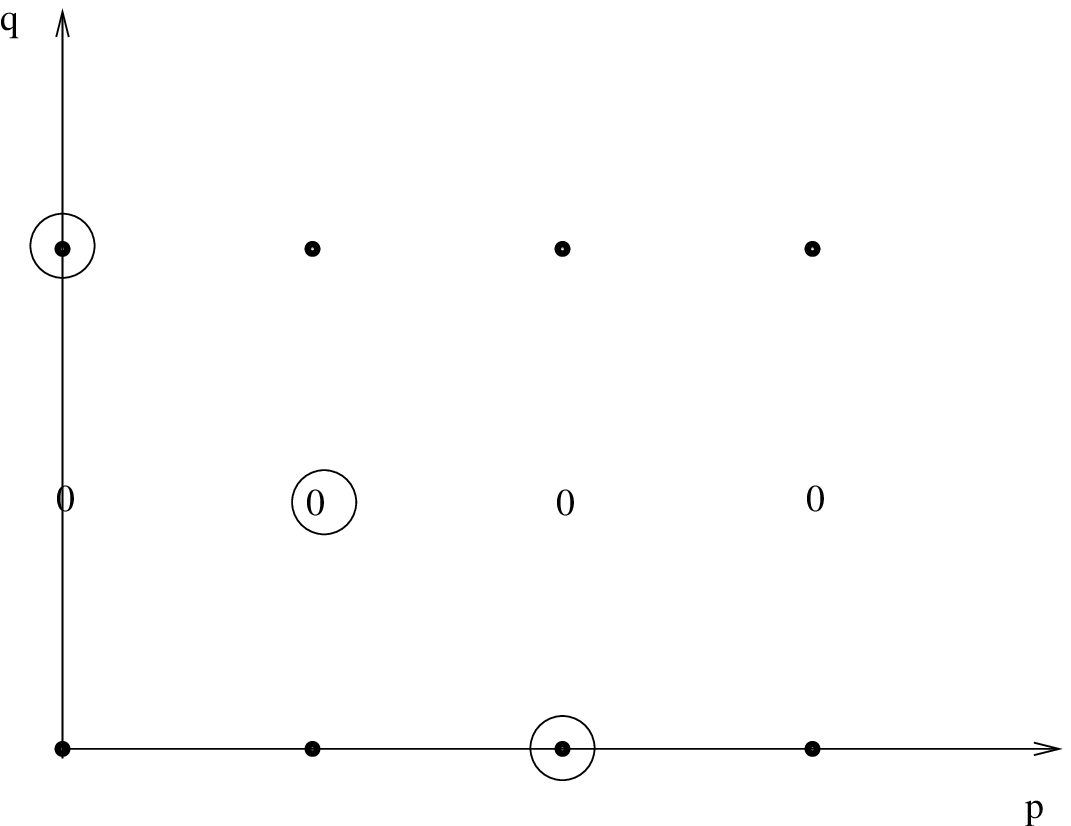}
  \caption{Terms $E^{pq}_{3}(\mathcal{X})$ and $E^{pq}_{4}(\mathcal{X})$}\label{fig:E3,E4}
  \end{center}
  \end{figure}
  
  In Figures \ref{fig:E1,E2} and \ref{fig:E3,E4}, the circled terms
  mean that they will stay constant that is they are equal to
  $E_{\infty}^{pq}(\mathcal{X})$.  We deduce that
  $E_{\infty}^{20}(\mathcal{X})=E_{2}^{20}(\mathcal{X})$ and
  $E_{\infty}^{02}(\mathcal{X})=E_{4}^{02}(\mathcal{X})$.

 The same argument for $T$ proves the lemma.
 \end{proof}

 \begin{proof}[Proof of Lemma \ref{lem:exact,seq,injective}]
 First, we show that the morphism $\alpha:E_{4}^{20}(\mathcal{X})\to
 E_{4}^{20}(T)$ is injective.
 From Figures \ref{fig:E1,E2} and \ref{fig:E3,E4}, we have that
 \begin{align}
   E_{4}^{20}(\mathcal{X})&=\ker\left(d_{3}:E_{2}^{20}(\mathcal{X})\to
   E_{2}^{03}(\mathcal{X})\right) \label{eq:2}\\
  E_{4}^{20}(T)&=\ker\left(d_{3}:E_{2}^{20}(T)\to
   E_{2}^{03}(T)\right)\label{eq:3}
 \end{align}
 Moreover, we have that
 \begin{align}
   E_{2}^{20}(\mathcal{X})&=\ker\left(d_{1}:H^{2}_{\et}(X_{0},\GG)\to
   H^{2}_{\et}(X_{1},\GG)\right)\label{eq:4}\\
  E_{2}^{20}(T)&=\ker\left(d_{1}:H^{2}_{\et}(T_{0},\GG)\to
   H^{2}_{\et}(T_{1},\GG)\right)\label{eq:5}
 \end{align}
 Recall that $U_{i}\simeq \aaa^{1}\times (\cc^{\ast})^{n-1}$ and
 $T_{0}=(\cc^{\ast})^{n}$. By Grothendieck's \textit{Expos\'es} on the Brauer group
 \cite[\S 6 p.133]{Ggb68}, we have the following long exact sequence :
 \begin{align}\label{eq:8}
  \xymatrix{ \ldots \ar[r] & H^{2}_{X_{0}-T_{0}}(X_{0},\GG)\ar[r]& H^{2}_{\et}(X_{0},\GG)\ar[r]& H^{2}_{\et}(T_{0},\GG) \ar[r] & \ldots }
 \end{align}
 Moreover, we have that
 \begin{itemize}
 \item the spectral sequence
 $F^{pq}_{2}:=H^{p}((X_{0}\setminus
 T_{0}),\underline{H^{q}_{(X_{0}\setminus T_{0})}(X_{0},\GG)})$ converges
 to $H_{(X_{0}\setminus T_{0})}(X_{0},\GG)$.
 \item  $H^{0}_{(X_{0}\setminus
   T_{0})}(X_{0},\GG)=H^{2}_{(X_{0}\setminus T_{0})}(X_{0},\GG)=0$
   and $H^{1}_{(X_{0}\setminus T_{0})}(X_{0},\GG)=\zz$.
 \end{itemize}
 This implies that $F_{2}^{20}=F_{2}^{02}=0$.
 As $X_{0}\setminus T_{0}=(\cc^{\ast})^{n-1}$, we have that
 $F_{2}^{11}=H^{1}(X_{0} \setminus T_{0},\zz)=0$.
 The spectral
 sequence $F_{2}^{pq}$ implies
 $H^{2}_{(X_{0}\setminus T_{0})}(X_{0},\GG)=0$. Hence, Sequence (\ref{eq:8}) and
 Equalities (\ref{eq:2}),~(\ref{eq:3}),~(\ref{eq:4}) and (\ref{eq:5}),
 imply that $\alpha$ is injective.

 Let prove that $\beta:E_{2}^{02}(\mathcal{X})\to
 E_{2}^{02}(T)$ is injective.  Recall that
 $E_{2}^{02}(\mathcal{X})=\ker d_{2}/\im d_{1}$ and
 $E_{2}^{02}(T)=\ker \widetilde{d}_{2}/\im \widetilde{d}_{1}$.
 We have the following commutative diagram
 \begin{displaymath}
   \xymatrix{H^{0}(X_{1},\GG)\ar[r]^{\delta_{1}}\ar@{^{(}->}[d]   &   H^{0}(X_{2},\GG)\ar@{^{(}->}[d]\ar[r]^{\delta_{2}}   &   H^{0}(X_{3},\GG)\ar@{^{(}->}[d] \\
 H^{0}(T_{1},\GG)\ar[r]^{\widetilde{\delta}_{1}}  &   H^{0}(T_{2},\GG)\ar[r]^{\widetilde{\delta}_{2}}  &   H^{0}(T_{3},\GG) }
 \end{displaymath}
 As $T_{ii'}\hookrightarrow U_{ii'}$ is open and dense, the vertical maps are
 injective. Notice that these maps are isomorphisms except on $U_{ii}$ and
 $U_{iii}$.  Let $\widetilde{y}\in H^{0}(T_{ii},\GG)$ such that
 there exists $x \in H^{0}(U_{iii},\GG)$ that lifts
 $\widetilde{\delta}_{1}(\widetilde{y})$ \ie we have the following diagram :
 \begin{displaymath}
   \xymatrix{& x \ar@{|->}[d]\\ \widetilde{y} \ar@{|->}[r] & \widetilde{\delta}_{1}(\widetilde{y})}
 \end{displaymath}
   The morphism $\widetilde{\delta}_{2}\vert_{T_{ii}} :
  H^{0}(T_{ii},\GG) \to H^{0}(T_{iii},\GG)$ is
  defined, for any $\widetilde{y}\in  H^{0}(T_{ii},\GG)$ and any
  $t,g,h \in T_{iii}=T_{i}\times \mu_{a_{i}}\times\mu_{a_{i}}$, by
  \begin{displaymath}
    \widetilde{\delta}_{2}\vert_{T_{ii}}(\widetilde{y})(t,g,h)=\widetilde{y}(ht,g)\widetilde{y}(t,h)/\widetilde{y}(t,gh).
  \end{displaymath}
 The divisor $U_{i}\setminus T_{i}$ is a principal divisor associate to
 the rational function $\varphi$.
 For any $g\in\mu_{a_{i}}$, the function $\widetilde{y}$ is rational
 on $U_{ii}\vert_{g}=U_{i}\times\{g\}$.
 Hence there exists a unique $n(g)$ in $\nn^{\ast}$ such that
 $\widetilde{y}\varphi^{n(g)}$ is a regular function on
 $U_{i}\times\{g\}$.
 As $\widetilde{y}(ht,g)\widetilde{y}(t,h)/\widetilde{y}(t,gh)$ is a
 regular function, we deduce that $\varphi^{n(g)+n(h)-n(gh)}=1$.
 Hence, the function $n:\mu_{a_{i}} \to \zz$ is a group homomorphism, therefore $n(g)=1$ for every $g$.
 We deduce that $\widetilde{y}$ is a regular function on $U_{i}$
 which implies that the morphism  $\beta:E_{2}^{02}(\mathcal{X})\to
 E_{2}^{02}(\mathcal{T})$ is injective.
 \end{proof}

\subsection{Characterization of a toric \DM stack as a gerbe over its rigidification}\label{sec:char-toric-dm}

Let $\XX$ be a toric \DM stack with \DM torus $\mathcal{T}$ isomorphic
to $T\times BG$ and coarse moduli space $X$. Denote by $\Xrig$ the
rigidification of $\XX$ (\cf Section \ref{subsec:rigidification})
which is by definition an orbifold with coarse moduli space $X$. The
universal property of the rigidification and of the canonical stack
(see Proposition \ref{prop;universal} and Corollary \ref{fmap}) imply
that we have the following strictly commutative diagram:
\begin{equation}\label{eq:18}
  \xymatrix{\XX \ar[r]^{r} \ar[d]_{f} & \Xrig  \ar[dl]^{\frig} \\ \XX^{\can} }
\end{equation}

Section \ref{subsec:rigidification} and  Lemma \ref{lem:igen,G,per,X}  imply that we can define
$\XX\fs G$. 

\begin{lem}\label{lem:orb,stru,rig} Let $\mathcal{X}$ be a toric \DM stack with \DM torus $\mathcal{T}$
  isomorphic to $T\times BG$.
   \begin{enumerate}
  \item\label{item:17} The orbifold $\Xrig$ is canonically isomorphic to $\XX\fs G$.
  \item\label{item:18}  There is a unique structure of toric orbifold on  $\Xrig$ with torus $T$
    such that the morphism $r:\XX \to \Xrig$ is a morphism of toric \DM
    stacks induced by $\TT\to T$.
  \end{enumerate}
\end{lem}

\begin{rem}\label{rem:lem,orb,stru,rig} Let $\mathcal{X}$ be a toric \DM stack with \DM torus
    $\mathcal{T}$ isomorphic to $T\times G$ and coarse moduli space $X$.
  \begin{enumerate}
  \item Proposition \ref{prop:morph,orb,can,toric} implies that the
    morphism $f^{\rig}:\Xrig\to \XX^{\can}$ is a morphism of toric \DM
    stacks. Hence we deduce that the commutative diagram (\ref{eq:18})
    is a commutative diagram of morphisms of toric \DM stacks.
\item  Let $H$ be a subgroup of $G$. The stack
  $\XX \fs H$ is a toric \DM stack with \DM torus isomorphic to
  $\mathcal{T} \fs H \simeq T\times \mathcal{B}(G/H)$. Moreover, the natural
  morphism $\XX \to \XX \fs H$ and $\XX \fs H \to \XX \fs G$ are
  morphism of toric \DM stacks.
\item Note that we did not use the non canonical isomorphism $\TT\equiv T\times\mathcal{B}G$ but only the short exact sequence of Picard stacks $1\to\mathcal{B}G\to\TT\to T\to 1$.
\end{enumerate}
\end{rem}

 \begin{proof}[Proof of Lemma \ref{lem:orb,stru,rig}] 
    (\ref{item:17}). As $\TT\fs G$ is isomorphic to the scheme
   $T$ which is open
   and dense in $\XX\fs G$, the stack $\XX\fs G$ is an orbifold which
   is canonically isomorphic to $\Xrig$.
 
    (\ref{item:18}).  The morphisms $\iota :\TT
   \hookrightarrow \XX$ and $a :\TT\times \XX \to \XX$ induce
   morphisms on the rigidifications $\iota^{\rig}:\TT\fs G \simeq T
   \to \Xrig$ and $ a^{\rig}:T \times \Xrig\to \Xrig$, by the
   universal property of the rigidification (See Proposition
   \ref{prop;universal}).  It is immediate to verify that $ a^{\rig}$
   is an action, extending the action of $T$ on itself. As $r^{-1}(T)$
   is isomorphic to $\TT$, we deduce that this is the only
   toric structure on $\Xrig$ which is compatible with the morphism $r$.
 \end{proof}

Since the morphism $r: \XX\to \Xrig$ is \'etale, the divisor
multiplicities of $\XX$ and $\Xrig$ are the same.

  \begin{thm}\label{thm:cara,gerbe,toric}
    \begin{enumerate}
    \item\label{item:20} Let $\YY$ be a toric orbifold with \DM torus
      $T$. Let $\XX\to \YY$ be an essentially trivial $G$-gerbe. Then
      $\XX$ has a unique structure of toric \DM stack with \DM torus
      isomorphic to $T\times \mathcal{B}G$ such that the morphism $\XX
      \to \YY$ is a morphism of toric \DM stacks.
    \item\label{item:19} Conversely, let $\XX$ be a toric \DM stack
      with \DM torus $\mathcal{T}\simeq T \times \mathcal{B}G$. Then
      $\XX\to \Xrig$ is an essentially trivial $G$-gerbe.
    \end{enumerate}
  \end{thm}

  \begin{proof}  (\ref{item:20}). The inverse
    image of $T$ in $\XX$, denoted by $\TT$, is open dense. 
    The restriction of the essentially trivial $G$-banded gerbe
    $\XX\to\YY$ to $T$ is the essentially trivial $G$-banded gerbe
    $\TT\to T$. Remark \ref{rem:gerbe}.(\ref{item:10}) implies that
    the gerbe $\TT\to T$ is trivial. The action of $T$ on $\YY$
    induces by pullback an action of $\TT$ on $\XX$. This is the only
    structure of toric \DM stack on $\XX$ compatible with the morphism
    $\XX\to \YY$.

     (\ref{item:19}).  Denote by $\alpha\in H^{2}_{\et}(\Xrig,G)$
    the $G$-banded gerbe $\XX\to \Xrig$. By Proposition
    \ref{prop:T,BG}, the restriction of $\alpha$ on the \DM torus $\TT$
    is the trivial $G$-banded gerbe in
    $H^{2}_{\et}(T,G)$.  Fix a cyclic decomposition of
    $G=\prod_{j=1}^{\ell}\mu_{b_{j}}$. By the isomorphism
    (\ref{eq:19}), the class $\alpha$ is sent to $(\alpha_{1}, \ldots
    ,\alpha_{\ell})\in \oplus_{j=1}^{\ell}
    H^{2}_{\et}(\Xrig,\mu_{b_{j}})$. We have that for any $j\in\{1,
    \ldots ,\ell\}$, the class of $\alpha_{j}$ restricts to the
    trivial class in $H^{2}_{\et}(T,\mu_{b_{j}})$. Theorem
    \ref{thm:gerbe,banale} states the injectivity of $\iota^{\ast}$
    in the following diagram.
\begin{displaymath}
  \xymatrix@1{
    H^1_{\et}(\Xrig,\mathbb{G}_m)
    \ar[rr]^-{\sqrt[b_j]{\cdot/\Xrig}}\ar[d] 
&& H^2_{\et}(\Xrig,\mu_{b_j}) \ar[d] \ar[rr] &&
 H_{\et}^2(\Xrig,\mathbb{G}_m)\ar@{^{(}->}[d]^{\iota^{\ast}} \\
  1 \ar[rr] && H_{\et}^2(T,\mu_{b_j}) \ar[rr] && H_{\et}^2(T,\mathbb{G}_m)
}
\end{displaymath}
A simple diagram chasing finishes the proof.
  \end{proof}


\begin{cor}\label{cor:carac,gerbe,toric}
  Let $\XX$ be a toric \DM stack with \DM torus $\TT$ isomorphic to $T\times
  \mathcal{B}G$.
  \begin{enumerate}
  \item 
Given $G=\prod_{j=1}^{\ell} \mu_{b_{j}}$. There exists
  ${L}_{j}$ in $\Pic(\Xrig)$ such that  $\mathcal{X}$ is isomorphic
  as $G$-banded gerbe over $\Xrig$ to 
  \begin{displaymath}
    \sqrt[b_{1}]{{L}_{1}/\Xrig}\times_{\Xrig}\cdots\times_{\Xrig}\sqrt[b_{\ell}]{{L}_{\ell}/\Xrig}.
  \end{displaymath}
  Moreover, the classes $([L_{1}], \ldots ,[L_{\ell}])$ in
  $\prod_{j=1}^{\ell}\Pic(\Xrig)/b_{j} \Pic(\Xrig)$ are unique.
\item The reduced closed substack $\mathcal{X}\setminus \mathcal{T}$ is a simple
normal crossing divisor.
  \end{enumerate}
\end{cor}

The first part of the corollary is very similar to Proposition 2.5 of \cite{Perroni-notetoricDeligne-Mumford-2007}.
\begin{rem}\label{rem:picard,group,gerbe} Let $\XX$ be a toric \DM
  stack with \DM torus $\TT$ isomorphic to $T\times \mathcal{B}G$.
  and $G=\prod_{j=1}^{\ell} \mu_{b_{j}}$.  Diagram (\ref{eq:21})
  and the corollary above imply that we have the following morphism of
  short exact sequences
\begin{equation}\label{eq:30}
\xymatrix{0\ar[r]&\zz^{\ell}\ar[r]^{\times \bs{b}}\ar[d]&\zz^{\ell}\ar[r]\ar[d]&\oplus_{j=1}^{\ell} \zz/b_{j}\zz\ar[r]\ar@{=}[d]&0\\0\ar[r]& \mathbb\Pic(\Xrig)\ar[r]^-{r^{\ast}} &\Pic(\mathcal{X})\ar[r]&\oplus_{j=1}^{\ell} \zz/b_{j}\zz\ar[r]&0}
\end{equation}
where the vertical morphisms sends $e_{j} \mapsto L_{j}$ and
$e_{j}\mapsto {L}^{1/b_{j}}_{j}$.
\end{rem}

\begin{proof}[Proof of Corollary \ref{cor:carac,gerbe,toric}]
 Theorem \ref{thm:cara,gerbe,toric}.(\ref{item:19}) implies that
 $\XX\to \Xrig$ is an essentially trivial $G$-banded gerbe. The first
 statement follows from Proposition \ref{prop:ess,trivi,gerb}.

By Corollary \ref{cor:reduced,toric}, we have that the reduced closed
substack $\Xrig\setminus \Trig$ is a simple normal crossing divisor. As the
morphism $\XX\to\Xrig$ is \'etale, we deduce the second statement of
the corollary.
\end{proof}

\section{Toric \DM stacks versus stacky fans}\label{sec:toric-dm-stacks-2}
In this Section, we will show that the toric \DM stacks that we have
defined correspond exactly with those of \cite{BCSocdms05}. 

In the first subsection, we show that our toric \DM stacks with a spanning
condition are global quotients. The second subsection makes the correspondence
with the article of \cite{BCSocdms05}.

\subsection{Toric \DM stacks as global quotients}\label{subsec:toric-dm-stacks-3}


Let $Z$ be a subvariety in $\cc^{n}$ of codimension
equal or higher than two. Let $G$ be
an abelian group scheme over $\cc$ that acts on $Z$ such that $[Z/G]$
is a \DM stack. According to Remark \ref{rem:sheaves,global,quotient},
a line bundle on $[Z/G]$ is given by a character $\chi$ of $G$.  Hence
the data of an invertible sheaf $L$ with a global section $s$ on
$[Z/G]$ give a morphism of groupoids between $[Z/G]$ and
$[\aaa^{1}/\cc^{\ast}]$.  Explicitly, this morphism is given by
$(s,\chi):Z\times G\to\aaa^{1}\times \cc^{\ast}$ and $s:Z\to
\aaa^{1}$.

In the following lemma, we use a slightly more general notion of a
root of  Cartier divisors that is a root of  invertible sheaves with 
global sections. All the properties of Section
\ref{subsubsec:root-line-section} are still true (see \cite{Cstc-2007}
or \cite{AGVgwdms}).

\begin{lem}\label{lem:root}
  Let $Z$ be a scheme. Let $G$ be an abelian
  group scheme over $\cc$ that acts on $Z$ such that $[Z/G]$ is a \DM
  stack. Let $(\bs{L},\bs{s}):=((L_{1},s_{1}), \ldots
  ,(L_{k},s_{k}))$ be $k$ invertible sheaves with global
  sections  over the quotient stack $[Z/G]$. Denote
by  $\bs{\chi}:=(\chi_{1}, \ldots ,\chi_{k})$ the representations
  associated to the invertible sheaves $\bs{L}$. Let $\bs{d}:=(d_{1},
  \ldots ,d_{k})$ be in $(\nn_{>0})^{k}$.
  \begin{enumerate}
  \item\label{item:4} We have that
    $\sqrt[\bs{d}]{(\bs{L},\bs{s})/[Z/G]}$ is isomorphic to $[\widetilde{Z}/\widetilde{G}]$
    where $\widetilde{Z}$ and $\widetilde{G}$ are defined by the following cartesian
    diagrams
    \begin{displaymath}
      \xymatrix{\widetilde{Z}\ar[r]\ar@{}[rd]|{\square}\ar[d]&\aaa^{k}
\ar[d]^-{\wedge \bs{d}}&&\widetilde{G}\ar[r]\ar[d]^-{\varphi}\ar@{}[rd]|{\square}
&\Gm^{k}\ar[d]^-{\wedge \bs{d}}\\
Z\ar[r]^-{\bs{s}}&\aaa^{k}&&G\ar[r]^-{\bs{\chi}}&\Gm^{k}}
    \end{displaymath} The action of $\widetilde{G}$ on $\widetilde{Z}$ is given by
    \begin{displaymath}  
(g,(\lambda_{1}, \ldots ,\lambda_{k})\cdot(z,(x_{1}, \ldots ,x_{k}))=(gz,(\lambda_{1} x_{1}, \ldots ,\lambda_{k}x_{k})
\end{displaymath}
for any $(g,(\lambda_{1}, \ldots ,\lambda_{k}))\in
\widetilde{G}$ and $(z,(x_{1}, \ldots
,x_{k}))\in \widetilde{Z}$.
\item\label{item:3} We have that $\sqrt[\bs{d}]{\bs{L}/[Z/G]}$ is
  isomorphic to $[Z/\widetilde{G}]$ where $\widetilde{G}$ is defined
  above. The action of $\widetilde{G}$ on $Z$ is given via $\varphi$.
  \end{enumerate}
\end{lem}

\begin{proof}
  It is a straightforward computation on fibered products of groupoids.
\end{proof}

\begin{rem}\label{rem:global,quotient}
  \begin{enumerate}
  \item \label{item:23} We have that $\ker \varphi$ is isomorphic to
    $\prod_{i=1}^{k}\mu_{d_{i}}$.  Notice that the action of
    $\widetilde{G}$ on $Z$ in the second part of the proposition above
    implies that the kernel of $\varphi$ acts trivially on $Z$. Hence,
    $[Z/\widetilde{G}]$ is a $\prod_{i=1}^{k}\mu_{d_{i}}$-banded gerbe
    over $[Z/G]$.
   \item\label{item:24} In both cases we have that
     $\widetilde{G}\in\Ext^{1}(G,\prod_{i=1}^{k}\mu_{d_{i}})$.
 \end{enumerate}
\end{rem}

\begin{lem}\label{lem:cocartesian}
  Let $A$ be an abelian group of finite type. Let $E$ in
  $\Ext^{1}(\oplus_{i=1}^{k}\zz/d_{i}\zz,A)$. If we have a morphism
  of short exact sequences:
  \begin{displaymath}
    \xymatrix{0\ar[r]& \zz^{k}\ar[rr]^{(d_{1}, \ldots ,d_{k})}\ar[d]&& \zz^{k}\ar[r]\ar[d]&
       \bigoplus_{i=1}^{k}\zz/d_{i}\zz \ar[r] \ar[d]^{\wr}&0\\
 0\ar[r]& A\ar[rr]& &E\ar[r]& \bigoplus_{i=1}^{k}\zz/d_{i}\zz \ar[r]&0}
  \end{displaymath} then the left square is cocartesian.
\end{lem}

\begin{rem} Diagrams (\ref{eq:28}) and (\ref{eq:30}) imply that we have the
       following cocartesian diagrams:
   \begin{equation}\label{eq:34}
     \xymatrix{\zz^{\ell}\ar[r]^{\times \bs{b}}\ar[d]& \zz^{\ell}
 \ar[d] & \zz^{n}\ar[r]^{\times \bs{a}}\ar[d]& \zz^{n} \ar[d]\\
 \Pic(\XX)\ar[r]& \Pic(\sqrt[\bs{b}]{\bs{L}/\XX})& \Pic(\XX)\ar[r]& \Pic(\sqrt[\bs{a}]{\bs{D}/\XX})}
\end{equation}
\end{rem}
 
\begin{proof}[Proof of Lemma \ref{lem:cocartesian}]
  Denote by $P$ the push-out of $\zz^{k}\to\zz^{k}$ and $\zz^{k}\to A$.
  Using the universal property of co-cartesian diagrams we deduce a
  morphism $f$ from $P$ to $E$ and the following morphisms of
  extensions:
\begin{displaymath}
  \xymatrix@1{
      0 \ar[r] & \zz^{k} \ar[rr]^-{\times(d_{1}, \ldots ,d_{k})}\ar[d] && \zz^{k} \ar[r]\ar[d] & \bigoplus_{i=1}^{k} \zz/{d_{i}}\zz \ar[r]\ar[d]^-{\alpha} & 0 \\
      0 \ar[r] & A\ar[rr]^{q}\ar@{=}[d] && P \ar[r]\ar[d]_-{f} & \coker(q)\ar[r]\ar[d]^-{\beta} & 0 \\
      0 \ar[r] & A \ar[rr]& & E \ar[r] & \bigoplus_{i=1}^{k} \zz/{d_{i}}\zz \ar[r] & 0 \\
}
\end{displaymath}
Notice that the composition $\beta\circ\alpha$ is the isomorphism in
Lemma \ref{lem:cocartesian}.
By simple diagram chasing, we deduce that $f$ is an isomorphism.
\end{proof}

\begin{rem}\label{rem:Z/G,\DM,toric}
  Let $\XX$ be a toric \DM stack with coarse moduli space
$X$. Proposition \ref{prop:toric,orbifold,cms,toric,variety} implies
that $X$ is a simplicial toric variety. Denote by $\Sigma$ a fan of $X$.
Assume that the rays of $\Sigma$ generate $N_{\qq}$.
As explained in Section \ref{subsec:geometry-toric}, we have that  $X$ is the
geometric quotient $Z_{\Sigma}/G_{A}$ where
$G_{A}:=\Hom(A^{1}(X),\cc^{\ast})$.
Put $G_{\XX}:=\Hom(\Pic(\XX),\cc^{\ast})$. Notice that $G_{\Xrig}$
acts on $Z_{\Sigma}$ via the dual (in the sense of Section
\ref{subsec:diag-group-scheme}) of the morphism $\zz^{n}\to
\Pic(\Xrig)$.  The group $G_{\XX}$ acts on $Z_{\Sigma}$ via the dual
of the morphism $\Pic(\Xrig) \to \Pic(\XX)$. Consider the quotient
stack $[Z_{\Sigma}/G_{\XX}]$. The quotient stack
$[(\ccs)^{n}/G_{\XX}]$ is a \DM torus which is open and dense in
$[Z_{\Sigma}/G_{\XX}]$. As the natural action of $(\ccs)^{n}$ on
$Z_{\Sigma}$ extends the action of $(\ccs)^{n}$ on itself, we deduce a
stack morphism $a: [(\ccs)^{n}/G_{\XX}]\times[Z_{\Sigma}/G_{\XX}] \to
[Z_{\Sigma}/G_{\XX}]$ that extends the action of
$[(\ccs)^{n}/G_{\XX}]$ on itself. Proposition
\ref{prop:action,morphism} implies that the stack morphism $a$ induces
a natural action of the \DM torus on $[Z_{\Sigma}/G_{\XX}]$ that is
$[Z_{\Sigma}/G_{\XX}]$ is a toric \DM stacks.
\end{rem}

\begin{thm}\label{thm:toric,global,quotient}
  Let $\XX$ be a toric \DM stack with coarse moduli space $X$. Denote
by  $\Sigma$ the fan associated to $X$. Assume that the rays of $\Sigma$
generate $N\otimes\qq$. Then $\XX$ is naturally isomorphic, as a toric
  stack, to $[Z_{\Sigma}/G_{\XX}]$ where $G_{\XX}:=\Hom(\Pic(\XX),\cc^{\ast})$. 
\end{thm}

\begin{rem}Removing the spanning condition of the rays gives the following
  result. Let $\XX$ be a toric \DM stack with torus $\TT$
  (isomorphic to $T\times \mathcal{B} G$) and with coarse moduli space the
  simplicial toric variety $X$. Denote by $\Sigma\subset N_{\qq}$ the fan of
  $X$.  From the footnote \ref{footnote} of Section \ref{subsec:geometry-toric},
  we deduce that the toric variety $X$ is isomorphic to $\widetilde{X}\times
  \widetilde{T}$ where $\widetilde{X}$ is a simplicial toric variety whose the rays
  of its fan $\widetilde{\Sigma}$
  span $\widetilde{N}_{\qq}$. Notice that the dimension of $\widetilde{T}$ is
  $\rank({N}_{\qq})-\rank(\widetilde{N}_{\qq})$. The previous Theorem implies
  that  $\XX$ is isomorphic, as toric stacks, to
  $[Z_{\widetilde{\Sigma}}/G_{\widetilde{\XX}}]\times (\widetilde{T}\times \mathcal{B}G)$
\end{rem}

  \begin{proof}[Proof of Theorem \ref{thm:toric,global,quotient}]If $\XX$ is $\XX^{\can}$, the theorem follows from Remark
    \ref{rem,can,pic}.(\ref{item:22}).  If $\XX$ is $\Xrig$, the
    theorem follows from the right cocartesian square of Diagram (\ref{eq:34})
    and Lemma \ref{lem:root}.(\ref{item:4}).  For a general $\XX$, it
    follows from the left cocartesian square of Diagram (\ref{eq:34}) and Lemma
    \ref{lem:root}.(\ref{item:3}).
  \end{proof}

\subsection{Toric \DM stacks and stacky fans}\label{sec:toric-dm-stacks}

First we recall the definition of a stacky fan from \cite{BCSocdms05}.
\begin{defn}\label{def:stacky,fan}
  \textit{A stacky fan} is a triple $\bs{\Sigma}:=(N,\Sigma,\beta)$
  where $N$ is a finitely generated abelian group, $\Sigma$ is a
  rational simplicial fan in $N_{\qq}:=N\otimes_{\zz}\qq$ with $n$
  rays, denoted by $\rho_{1}, \ldots ,\rho_{n}$, and a morphism of
  groups $\beta: \zz^{n}\to N$ such that
  \begin{enumerate}
  \item the rays span $N_{\qq}$,
  \item for any $i\in\{1, \ldots ,n\}$, the element
    $\overline{\beta(e_{i})}$ in $N_{\qq}$ is on the ray $\rho_{i}$
    where $(e_{1}, \ldots ,e_{n})$ is the canonical basis of $\zz^{n}$
    and the natural map $N\to N_{\qq}$ sends $m\mapsto\overline{m}$.
 \end{enumerate}
\end{defn}

\begin{rem}\label{rem:stacky,fan} 
  Let $\bs{\Sigma}:=(N,\Sigma,\beta)$ be a stacky fan.
  \begin{enumerate}
  \item As the rays span $N_{\qq}$, we have that $\beta$ has finite
    cokernel.
  \item For any $i\in\{1, \ldots ,n\}$, denote by $v_{i}$ the unique
    generator of $\rho_{i}\cap (N/N_{\tor})$ where $N_{\tor}$ is the
    torsion part part of $N$. Denote by $\beta^{\rig}$ the composition of
    $\beta$ followed by the quotient morphism $N\to N/N_{\tor}$. There
    exists a unique $a_{i}\in \nn_{>0}$ such that
    $\beta^{\rig}(e_{i})=a_{i}v_{i}$. Denote 
    $\bs{\Sigma}^{\rig}:=(N/N_{\tor},\Sigma,\beta^{\rig})$. There
    exists a unique group homomorphism $\beta^{\can}:\zz^{n}\to
    N/N_{\tor}$ such that we have the following commutative diagram.
  \begin{equation}\label{eq:11}
    \xymatrix{\zz^{n} \ar[rr]^-{\beta} \ar[rrdd]^{\beta^{\rig}} \ar[dd]_-{\diag(a_{1}, \ldots ,a_{n})}&& N \ar[dd] \\ \\ \zz^{n} \ar[rr]^-{\beta^{\can}} && N/N_{\tor}}
  \end{equation}
  Denote $\bs{\Sigma}^{\can}:=(N/N_{\tor},\Sigma,\beta^{\can})$.
  
  In Remark 4.5 of \cite{BCSocdms05}, the authors define the notion of
  morphism of stacky fans. The commutative diagram (\ref{eq:11})
  provides us the morphisms of stacky fans
  $\bs{\Sigma}\to\bs{\Sigma}^{\rig}\to\bs{\Sigma}^{\can}$.
\item To the fan $\Sigma$, we can associate canonically the stacky fan
  $\bs{\Sigma}^{\can}$.
  \end{enumerate}
\end{rem}

\begin{construction}[Construction of the \DM stack associated to the stacky fan
  $\bs{\Sigma}$]\label{const:stack,toric}
Now we explain how to associate a \DM stack $\XX(\bs{\Sigma})$ to a
stacky fan $\bs{\Sigma}$ following Sections $2$ and $3$ in \cite{BCSocdms05}.
Denote by $d$ the rank of $N$. Choose a projective resolution of
$N$ with two terms
that is 
\begin{displaymath}
  \xymatrix{0 \ar[r]& \zz^{\ell} \ar[r]^-{Q}& \zz^{d+\ell} \ar[r]& N \ar[r]&0}
\end{displaymath}
Choose a map $B:\zz^{n}\to
\zz^{d+\ell}$ lifting the map $\beta:\zz^{n}\to N$.
Consider the morphism $[BQ]:\zz^{n+\ell}\to \zz^{d+\ell}$.
Denote $DG(\beta):=\coker([BQ]^{\ast})$. Denote by 
$\beta^{\vee}:(\zz^{n})^{\ast}\to DG(\beta)$ the group morphism that
makes the following diagram commute
\begin{displaymath}
  \xymatrix{(\zz^{n})^{\ast} \ar@{^{(}->}[r] \ar[rd]_-{\beta^{\vee}}&
    (\zz^{n+\ell})^{\ast}\ar[d]
 \\ & DG(\beta):=\coker{[BQ]^{\ast}}}
\end{displaymath}

Let $Z_{\Sigma}$ be the quasi-affine variety  associated to the
fan $\Sigma$ (see Section \ref{subsec:geometry-toric}). Define the
action of $G_{\bs{\Sigma}}:=\Hom_{\zz}(DG(\beta),\cc^{\ast})$ on
$Z_{\Sigma}$ as follows.
Applying the functor $\Hom_{\zz}(-,\cc^{\ast})$ to the morphism
$\beta^{\vee}:(\zz^{n})^{\ast}\to DG(\beta)$, we get a group
morphism $G_{\bs{\Sigma}}\to(\cc^{\ast})^{n}$. Via the natural action of
$(\cc^{\ast})^{n}$ on $\cc^{n}$, we define an action of $G_{\bs{\Sigma}}$ on
$Z_{\Sigma}$. Finally, the stack associated to the stacky fan
$\bs{\Sigma}:=(N,\Sigma,\beta)$ is the quotient stack
$\mathcal{X}(\bs{\Sigma}):=[Z_{\Sigma}/G_{\bs{\Sigma}}]$.
\end{construction}

\begin{notn}
  We will later see that the group $G_{\bs{\Sigma}}$ is isomorphic to
   $G_{\XX}:=\Hom(\Pic(\XX),\cc^{\ast})$.
\end{notn}

 By Proposition 3.2 of \cite{BCSocdms05}, we have
that $[Z_{\Sigma}/G_{\bs{\Sigma}}]$ is a smooth \DM stack.

\begin{rem}
  In  \cite{Iwanari-Logarithmic-2007}, Iwanari defined a
  smooth toric Artin stack over any scheme associated to a stacky fan
  $\bs{\Sigma}^{\rig}$.
\end{rem}

\begin{rem}\label{rem:span,problem}
  As it was observed in Section 4 of \cite{BCSocdms05}, the condition
  that the rays span $N_{\qq}$ in Definition \ref{def:stacky,fan} is
  not natural. Indeed a \DM torus $(\ccs)^{d}\times \mathcal{B}G$ where
  $G$ is a finite abelian group can not be produced as a stack
  $\XX(\bs{\Sigma})$ for $\bs{\Sigma}$ a stacky fan with the condition
  that the rays span $N_{\qq}$. 
  Nevertheless, it is not really true to say that toric \DM stacks are a
  ``generalization" of the stacks $\XX(\bs{\Sigma})$. Indeed, as for toric
  variety, we will see that a toric \DM stack is a product of a
  $\XX(\bs{\Sigma})$ by a \DM torus.
\end{rem}

\begin{lem}\label{lem:stacky,fan}
  Let $\bs{\Sigma}:=(N,\Sigma,\beta)$ be a stacky fan.
  \begin{enumerate}
  \item\label{item:25} The stack $\XX(\bs{\Sigma})$ is a toric \DM stack.
  \item\label{item:26} The stack $\XX(\bs{\Sigma})$ is a toric orbifold if and only
  if the finitely generated abelian group $N$ is free. 
\item\label{item:27} The stack $\XX({\bs{\Sigma}})$ is canonical if and only if $\bs{\Sigma}=\bs{\Sigma}^{\can}$.
\end{enumerate}
\end{lem}

\begin{proof}(\ref{item:25}). The group morphism $G_{\bs{\Sigma}}\to
  (\cc^{\ast})^{n}$ defined in Construction \ref{const:stack,toric}
  defines the quotient stack $[(\cc^{\ast})^{n}/G_{\bs{\Sigma}}]$
  which is by definition a \DM torus. As the open dense immersion
  $(\cc^{\ast})^{n}\hookrightarrow Z_{\Sigma}$ is
  $G_{\bs{\Sigma}}$-equivariant, we have that the stack morphism
  $[(\cc^{\ast})^{n}/G_{\bs{\Sigma}}]\to [Z_{\Sigma}/G_{\bs{\Sigma}}]$
  is an open dense immersion. Using the same arguments of Remark
  \ref{rem:Z/G,\DM,toric}, we have that the action of the \DM torus
  $[(\cc^{\ast})^{n}/G_{\bs{\Sigma}}]$ on itself extends to an action
  on $[Z_{\Sigma}/G_{\bs{\Sigma}}]$. That is $\XX(\bs{\Sigma})$ is a
  toric \DM stack. 
  
(\ref{item:26}). The stack $\XX(\bs{\Sigma})$ is a toric
  orbifold if and only if $G_{\bs{\Sigma}}\to (\cc^{\ast})^{n}$ is
  injective, if and only if $\beta^{\vee}$ is surjective, if and only if
  $N$ is free.
  
(\ref{item:27}). Assume that $\bs{\Sigma}=\bs{\Sigma}^{\can}$.
  As the coarse moduli space $X$ of $\XX({\bs{\Sigma}})$ is the
  geometrical quotient $Z_{\Sigma}/G_{A^{1}(X)}$ where
  $G_{A^{1}(X)}:=\Hom(A^{1}(X),\cc^{\ast})$, we have that
  $\XX^{\can}=[Z_{\Sigma}/G_{A^{1}(X)}]$. Construction
    \ref{const:stack,toric} implies that $G_{\bs{\Sigma}}$ is
    $G_{A^{1}(X)}$.  Conversely, if $\bs{\Sigma}\neq\bs{\Sigma}^{\can}$
    then either $N$ has torsion (\ie $\XX(\bs{\Sigma})$ is a gerbe) or there exists a
    divisor $D$ associated to a ray such that any geometric point of
    $D$ has a non trivial stabilizer.
\end{proof}

\begin{rem}\label{rem:BCS,pic}
  Let $\XX(\bs{\Sigma})$ be a canonical stack (\ie
  $\bs{\Sigma}=\bs{\Sigma}^{\can}$). The proof of the third statement of Lemma
  \ref{lem:stacky,fan} implies that $DG(\beta^{\can})=\Pic(\XX(\bs{\Sigma}))$.
\end{rem}

\begin{thm}\label{thm:orb,stacky}
  Let $\XX$ be a toric orbifold with coarse moduli space $X$. Denote by
  $\Sigma$ a fan of $X$ in $N_{\qq}:=N\otimes_{\zz}\qq$. Assume that
  the rays of $\Sigma$ span $N_{\qq}$. Then there is a
  unique $\beta:\zz^{n}\to N$ such that the stack associated to the
  stacky fan $(N,\Sigma,\beta)$ is isomorphic as toric orbifold to
  $\XX$.
\end{thm}
                                
\begin{rem} An arbitrary toric orbifold is isomorphic to a product
$\XX(\bs{\Sigma})\times (\ccs)^{k}$. 
\end{rem}
                                
\begin{proof}[Proof of Theorem \ref{thm:orb,stacky}]
 Denote by $\bs{a}:=(a_{1}, \ldots ,a_{n})$
  the divisor multiplicities of $\XX$. We define the morphism of
  groups $\beta:\zz^{n}\to N$ by sending $e_{i}\mapsto a_{i}v_{i}$
  where $v_{i}$ is the generator of the semi-group $\rho_{i}\cap N$.
  Denote by  $\bs{\Sigma}$ the stacky fan $(N,\Sigma,\beta)$.  
  
  Theorem \ref{thm:toric,global,quotient} states that $\XX$ is
  isomorphic to $[Z_{\Sigma}/G_{\XX}]$. In order to prove that the two
  stacks are isomorphic, we will show that $G_{\XX}$ is
  isomorphic to $G_{\bs{\Sigma}}$ such that the two actions on
  $Z_{\Sigma}$ are compatible.
  From Diagram (\ref{eq:11}), we deduce a morphism of exact sequences:
  \begin{displaymath}
    \xymatrix{
0\ar[r]&(\zz^{n})^{\ast}\ar[r]\ar[d]&(\zz^{n})^{\ast}\ar[d]\ar[r]&\oplus_{i=1}^{n} \zz/a_{i}\zz\ar[r]\ar@{=}[d]&0\\
0\ar[r]&\Pic(\XX^{\can})\ar[r]&DG(\beta)\ar[r]&\ar[r]\oplus_{i=1}^{n} \zz/a_{i}\zz&0}
  \end{displaymath}
 The right cocartesian square of Diagram (\ref{eq:34}) implies that
  $G_{\bs{\Sigma}}$ is isomorphic to $G_{\XX}$ such
  that the actions of $G_{\bs{\Sigma}}$ and $G_{\XX}$ on
  $Z_{\Sigma}$ are compatible.
  
  The uniqueness of $\beta$ follows from the geometrical
   interpretation of the divisor multiplicities.
\end{proof}

  
  
  

\begin{rem}\label{rem:Pic,orb}
  \begin{enumerate}
  \item The proof shows also that $\Pic(\XX)$ is isomorphic to $DG(\beta)$.
  \item Marking a point $a_{i}v_{i}$ on the ray $\rho_{i}\cap N$
  corresponds geometrically to putting a generic stabilizer $\mu_{a_{i}}$
  on the divisor $\DD_{i}$ associated to the ray $\rho_{i}$.
  \end{enumerate}
\end{rem}

\begin{prop}\label{prop:stacky,rig}
Let $\bs{\Sigma}:=(N,\Sigma,\beta)$ be a stacky fan. There is a
unique $\alpha$ in \\
$\Ext^{1}(N_{\tor},\Pic(\XX(\bs{\Sigma}^{\rig})))$ such that the
essentially trivial $\Hom(N_{\tor},\cc^{\ast})$-banded gerbe over $\XX(\bs{\Sigma}^{\rig})$
associated to $\alpha$ is isomorphic as banded gerbe to $\XX(\bs{\Sigma})$.
\end{prop}

\begin{proof}Fix a decomposition $N=\zz^{d}\oplus \oplus_{j=1}^{\ell}
  \zz/b_{j}\zz$. It follows from Construction \ref{const:stack,toric}
  that we have the following  diagram:
  \begin{equation}\label{eq:25}
    \xymatrix{
&0&0&0&\\
0\ar[r]&DG(\beta^{\rig})\ar[r]\ar[u]&DG(\beta)\ar[r]\ar[u]& \oplus_{j=1}^{\ell}
  \zz/b_{j}\zz\ar[r]\ar[u]&0\\
0\ar[r]&\ar[r](\zz^{n})^{\ast}\ar[u]^{(\beta^{\rig})^{\vee}}&\ar[r]\ar[u](\zz^{n+\ell})^{\ast}&\ar[u]\ar[r](\zz^{\ell})^{\ast}&0\\
0\ar[r]&\ar[u]\ar[r](\zz^{d})^{\ast}&\ar[u]^{[BQ]^{\ast}}\ar[r](\zz^{d+\ell})^{\ast}&\ar[u]_{\times
(b_{1}, \ldots ,b_{\ell})}\ar[r](\zz^{\ell})^{\ast}&0\\
&0\ar[u]&\ar[u]0&\ar[u]0&}
\end{equation}
From Remark \ref{rem:Pic,orb}, we have that
$\Pic(X(\bs{\Sigma}^{\rig}))$ is isomorphic to $DG(\beta^{\rig})$.
The first line of Diagram (\ref{eq:25}) is an element $\alpha\in
\Ext^{1}(N_{\tor},\Pic(\XX(\bs{\Sigma}^{\rig})))$.  By Proposition
\ref{prop:ess,trivi,gerb}, the element $\alpha$ induces an element
$([L_{1}], \ldots ,[L_{\ell}])\in
\prod_{j=1}^{\ell}\Pic(\Xrig)/b_{j}\Pic(\Xrig)$.
 The last row of the diagram above is a projective resolution of
 $\oplus_{j=1}^{\ell}\zz/b_{j}$. Hence, we deduce that there exists a  morphism
 of short exact sequence 
 \begin{equation}\label{eq:32}
     \xymatrix{
 0\ar[r]&(\zz^{\ell})^{\ast}\ar[r]^{\times \bs{b}} \ar[d]^{\widetilde{f}}&(\zz^{\ell})^{\ast}\ar[d]^{f}\ar[r]&
 \oplus_{j=1}^{\ell} \zz/b_{j}\zz\ar[r]\ar@{=}[d]&0\\
 0\ar[r]&\Pic(\Xrig)\ar[r]&DG(\beta)\ar[r]&
 \ar[r]\oplus_{j=1}^{\ell} \zz/b_{j}\zz&0}
 \end{equation}
 The morphism $\widetilde{f}$ is the same as the choice of $L_{1},
 \ldots ,L_{\ell}$ in $\Pic(\Xrig)$ in the classes $[L_{1}], \ldots
 ,[L_{\ell}]$.
By the left cocartesian square of Diagram (\ref{eq:34}), we deduce
 that $G_{\bs{\Sigma}}$ is isomorphic to $G_{\XX}$.
 We conclude that $\XX$ is isomorphic to
  $\XX(\bs{\Sigma})$. The uniqueness of $\alpha$ follows from Proposition \ref{prop:ess,trivi,gerb}.


\end{proof}

\begin{rem}
  Denote  by $\XX_{1}$ and by $\XX_{2}$ respectively the stacks
  associated to stacky fans $(\Sigma,N,\beta_{1})$ and
  $(\Sigma,N,\beta_{2})$.  The stacks $\XX_{1}$ and $\XX_{2}$ are
  isomorphic, as toric \DM stack, if and only if the extensions defined in Diagram
  (\ref{eq:25}) in $\Ext^{1}(N_{\tor},\Pic(\XX(\bs{\Sigma}^{\rig})))$
  are isomorphic.
\end{rem}

\begin{thm}\label{thm:toric,BCS}
  Let $\XX$ be a toric \DM stack with coarse moduli space $X$. Denote
  by $\Sigma$ a fan of $X$ in $N_{\qq}$. Assume that the rays of $\Sigma$
  span $N_{\qq}$. There exist $N$ and $\beta :\zz^{n}\to N$ such that
  the stack associated to the stacky fan $(N,\Sigma,\beta)$ is
  isomorphic as toric \DM stacks to $\XX$.
  \end{thm}

  \begin{rem}
Let $\bs{\Sigma}$ be a stacky fan. Corollary \ref{cor:carac,gerbe,toric} 
and the theorem above imply that $\XX(\bs{\Sigma})$ is
 isomorphic to a product of root stacks over its rigidification. 
This result was discovered independently by Perroni 
(\cf Proposition 3.2 in \cite{Perroni-notetoricDeligne-Mumford-2007})
 and by Jiang and Tseng (\cf Remark 2.10 in \cite{Jiang-Tseng-IntegralChowRing-2007}).       
\end{rem}

 \begin{proof}[Proof of Theorem \ref{thm:toric,BCS}]
    If $\XX$ is a toric orbifold then the statement was already proved
    in Theorem \ref{thm:orb,stacky}.
    
    Let $\XX$ be a toric \DM stack with \DM torus isomorphic to $T\times
    \mathcal{B}G$. By Theorem \ref{thm:toric,global,quotient}, we have
    that $ \XX$ is isomorphic to $[Z_{\Sigma}/G_{\XX}]$. By Theorem
    \ref{thm:orb,stacky}, there exists a unique stacky fan
    $\bs{\Sigma}^{\rig}=(\Sigma,\zz^{d},\beta^{\rig})$ where $d:=\dim
    \XX$ such that $\Xrig$ is isomorphic to $\XX(\bs{\Sigma}^{\rig})$.
    
    There exist $(b_{1}, \ldots ,b_{\ell})\in (\nn_{>0})^{\ell}$ such that
    $G=\prod_{j=1}^{\ell}\mu_{b_{j}}$. Put $N:=\zz^{d}\oplus
    \oplus_{j=1}^{\ell}\zz/b_{j}\zz$. Corollary
    \ref{cor:carac,gerbe,toric} gives us $\ell$ invertible sheaves
    $L_{1}, \ldots ,L_{\ell}$ on $\Xrig$. For any $j$, choose 
    $c_{1j}, \ldots ,c_{nj} \in \zz$ such that
    $L_{j}=\otimes_{i=1}^{n} \mathcal{O}(\DD_{i}^{\rig})^{c_{ij}}$
    where $\DD_{i}^{\rig}$ is the Cartier divisor associated to the
    ray $\rho_{i}$. Put 
    \begin{align*}
      \beta : \zz^{n} &\longrightarrow \zz^{d}\oplus
    \oplus_{j=1}^{\ell}\zz/b_{j}\zz \\
      e_{i} &  \longmapsto (\beta^{\rig}(e_{i}),[c_{i1}], \ldots ,[c_{i\ell}])
    \end{align*} where $[c_{ij}]$ is the class of $c_{ij}$ modulo
    $b_{j}$. It is straightforward to check that $X(\bs{\Sigma})$ is
    isomorphic to $\XX$.
    \end{proof}

    \begin{rem}\label{rem:toric,choices} Let $\XX$ be a toric \DM stack with \DM torus isomorphic to
      $T\times \mathcal{B}G$.  The non-uniqueness of $N$ and $\beta$ comes from
      three different kinds :
      \begin{enumerate}
      \item the decomposition of $G$ in product of cyclic groups \ie
        $G=\prod_{j=1}^{\ell}\mu_{b_{j}}$,
      \item the choice of the lift for the class $[L_{j}]\in
        \Pic(\XX^{\rig})/b_{j}\Pic(\XX^{\rig})$ for $j=1, \ldots ,\ell$ and
      \item the choice of the decomposition $L_{j}=\otimes_{i=1}^{n}
        \mathcal{O}(\DD_{i}^{\rig})^{c_{ij}}$ (see Example \ref{expl:P(4,6)} for
        such an example).
      \end{enumerate}
    \end{rem}

\subsection{Examples}
\label{sec:examples}

\begin{exmp}[Weighted projective spaces]
  Let $w_{0}, \ldots ,w_{n}$ be in $\nn_{>0}$. Denote by $\pp(w)$ the quotient
  stack $[\cc^{n+1}\setminus \{0\}/\cc^{\ast}]$ where the action of $\cc^{\ast}$
  is defined by $\lambda(x_{0}, \ldots ,x_{n})=(\lambda^{w_{0}}x_{0}, \ldots
  ,\lambda^{w_{n}}x_{n})$ for any $\lambda\in \cc^{\ast}$ and any $(x_{0},
  \ldots ,x_{n})\in \cc^{n+1}\setminus\{0\}$. The stack $\pp(w)$ is a complete
  toric \DM stack with \DM torus $[(\cc^{\ast})^{n+1}/\cc^{\ast}]\simeq
  \cc^{n}\times B\mu_{d}$  where $d:=\gcd(w_{0}, \ldots
  , \ldots ,w_{n})$ (cf Example \ref{exmp:P(w),torus}).

We have that
\begin{enumerate}
\item the stack $\pp(w)$ is canonical if and only if for any $i\in\{0, \ldots
  ,n\}$, we have that $\gcd(w_{0}, \ldots ,\hat{w}_{i}, \ldots
  ,w_{n})=1$ (\eg the weights are well-formed).
\item The stack $\pp(w)$ is an orbifold if and only if $\gcd(w_{0}, \ldots
  , \ldots ,w_{n})=1$.
\item The Picard group of $\pp(w)$ is cyclic. More precisely, we have 
  \begin{displaymath}
    \Pic(\pp(w))=\begin{cases}
      \zz&\mbox{ if } \dim \pp(w)\geq 1\\
    \zz/w_{0}\zz&\mbox{ if } \pp(w)=\pp(w_{0}). 
    \end{cases}
  \end{displaymath}
\end{enumerate}

  \begin{prop}\label{prop:charac,wps}
    Let $\XX$ be a complete toric \DM stack of dimension $n$ such that
    its Picard group is cyclic.  Then there exists unique up to order $(w_{0}, \ldots ,w_{n})$ in
    $(\nn_{>0})^{n+1}$ such that $\XX$ is isomorphic to $\pp(w_{0},
    \ldots ,w_{n})$.
  \end{prop}

  \begin{proof}Denote by $X$ the coarse moduli space of $\XX$. Denote
by    $\Sigma$ a fan of $X$. If the Picard group is isomorphic to $\zz/d\zz$
    then Theorem \ref{thm:toric,global,quotient} implies that
    $\XX=[Z_{\Sigma}/\mu_{d}]$ with $Z_{\Sigma}\subset \cc^{n}$.
    Hence, the fan $\Sigma$ has $n$ rays. In this case, $\XX$ is complete if
    and only if $n=0$. We deduce that $\XX=\mathcal{B}\mu_{d}\simeq\pp(d)$.
    
    If the Picard group is $\zz$, Theorem
    \ref{thm:toric,global,quotient} implies that
    $\XX=[Z_{\Sigma}/\cc^{\ast}]$ with $Z_{\Sigma}\subset \cc^{n+1}$.
    As $X$ is complete, the fan $\Sigma$ is complete. We deduce that
    $Z_{\Sigma}=\cc^{n+1}\setminus\{0\}$. The \DM torus is isomorphic to
    $[(\cc^{\ast})^{n+1}/\cc^{\ast}]$. The action of $\cc^{\ast}$ is
    given by the morphism $\cc^{\ast}\to (\cc^{\ast})^{n+1}$ that
    sends $\lambda\mapsto (\lambda^{w_{0}}, \ldots ,\lambda^{w_{n}})$
    with $w_{i}\in \zz \setminus\{0\}$. Notice that if the $w_{i}$'s do not
    have the same sign then $\XX$ is not separated. If the $w_{i}$'s
    are all negative then replacing $\lambda$ by $\lambda^{-1}$
    induces an isomorphism with a weighted projective space.
  \end{proof}
\end{exmp}

\begin{exmp}\label{expl:P(4,6)}
In this example, we give two isomorphic stacky fans for
$\pp(6,4)$ which was considered in Example 3.5 in \cite{BCSocdms05}.
As we have seen in Section \ref{sec:toric-dm-stacks}, $N$ and $\Sigma$
are fixed whereas $\beta$ is not unique.
Let $N$ be $\zz\times\zz/2$. Let $\Sigma$ be the fan in $N_{\qq}=\qq$
where the cones are $0,\qq_{\geq 0},\qq_{\leq 0}$.
Put 
\begin{align}\label{eq:33}
  \beta_{1}:\zz^{2}&\longrightarrow \zz\times \zz/2 &  \beta_{2}:\zz^{2}&\longrightarrow \zz\times \zz/2 \\
  e_{1}&\longmapsto (2,1) &  e_{1}&\longmapsto (2,1)\nonumber \\
  e_{2}&\longmapsto (-3,0) &  e_{2}&\longmapsto (-3,1) \nonumber
\end{align}
One can check  that the stack associated to
$(N,\Sigma,\beta_{1})$ and $(N,\Sigma,\beta_{2})$ is $\pp(6,4)$.

Let us explicit the bottom up construction in this case.  Its coarse moduli
space is $\pp^{1}$. The rigidification of $\pp(6,4)$ is $\pp(3,2)$.  Denote by
$x_{1},x_{2}$ the homogeneous coordinates of $\pp^{1}$. We have that
$\pp(3,2)=\sqrt[(2,3)]{(D_{1},D_{2})/\pp^{1}}$ where $D_{i}$ is the Cartier
divisor $(\mathcal{O}_{\pp^{1}}(1),x_{i})$. We have that
$\mathcal{O}_{\pp(3,2)}(\DD_{1})=\mathcal{O}_{\pp(3,2)}(3)$,
$\mathcal{O}_{\pp(3,2)}(\DD_{2})=\mathcal{O}_{\pp(3,2)}(2)$ and
$\pi^{\ast}\mathcal{O}_{\pp^{1}}(1)=\mathcal{O}_{\pp(3,2)}(6)$ where
$\pi:\pp(3,2) \to \pp^{1}$ is the structure morphism.  The stack $r:\pp(6,4)\to
\pp(3,2)$ is a $\mu_{2}$-banded gerbe isomorphic to
$\sqrt[2]{\mathcal{O}_{\pp(3,2)}(1)/\pp(3,2)}$.  In $\Pic(\pp(3,2))/2
\Pic(\pp(3,2))$, the class of $\mathcal{O}_{\pp(3,2)}(1)$ is also the class of
$\mathcal{O}_{\pp(3,2)}(\DD_{1})$ or the class of
$\mathcal{O}_{\pp(3,2)}(\DD_{1})\otimes\mathcal{O}_{\pp(3,2)}(\DD_{2})$.  These
different choices lead to the two isomorphic stacky fans in (\ref{eq:33}).
\end{exmp}

\begin{exmp}[Complete toric lines]\label{exmp:toric,line} Here, we explicitly describe all complete
  toric orbifolds of dimension $1$.  Notice that the coarse moduli
  space of a complete toric line is $\pp^{1}$. Denote by $x_{1},x_{2}$
  the homogeneous coordinates. Let $D_{i}$ the Cartier divisor
  $(\mathcal{O}(1),x_{i})$. Let $a_{1},a_{2}$ in $\nn_{>0}$. Denote by
  $d$ (resp. $m$) the greatest common divisor (resp. the lowest common
  multiple) of $a_{1},a_{2}$. The Picard group of the root stack
  $\sqrt[(a_{1},a_{2})]{(D_{1},D_{2})/\pp^{1}}$ is isomorphic to
  $\zz\times (\zz/d\zz)$. Notice that it is not a weighted projective
  space in general. As a global quotient, the stack
  $\sqrt[(a_{1},a_{2})]{(D_{1},D_{2})/\pp^{1}}$ is
  $[(\cc^{2}\setminus\{0\})/(\ccs\times\mu_{d})]$ where the action is
  given by
  \begin{align*}
    \ccs\times \mu_{d} \times(\cc^{2}\setminus\{0\}) &\longrightarrow
    (\cc^{2}\setminus\{0\}) \\
((\lambda,t),(x_{1},x_{2}))& \longmapsto (\lambda^{m/a_{1}}t^{k_{2}}x_{1},\lambda^{m/a_{2}}t^{-k_{1}}x_{2}) 
  \end{align*} 
  where $k_{1}, k_{2}$ are integers such that
  $\frac{k_{1}}{a}+\frac{k_{2}}{b}=\frac{1}{m}$. 
\end{exmp}

\appendix
\section{Uniqueness of morphisms to separated stacks}\label{sec:appendix,proof}
We prove proposition \ref{prop:sep,bis}.
\begin{prop}
  Let $\XX$ and $\YY$ be two \DM stacks. Assume that $\XX$ is normal
  and $\YY$ is separated.  Let $\iota:\mathcal{U}\hookrightarrow \XX$
  be a dominant open immersion. If $F_{1},F_{2}
  :\XX \to \YY$ are two morphisms of stacks such that there exits a
  $2$-arrow $\beta:F_{1}\circ \iota {\Rightarrow}F_{2}\circ \iota$
  then there exists a unique $2$-arrow $\alpha:F_{1}\Rightarrow F_{2}$ such
  that $\alpha\ast \Id_{\iota}=\beta$.
\end{prop}

\begin{proof}
  \textbf{Uniqueness:} We first assume that $\XX$ is a scheme, denoted
  by $X$, and $\YY$ is a global quotient $[V/G]$ where $G$ is a
  separated group scheme. Denote by $U$ the scheme $\UU$, open dense in
  $X$. For $i$ in $\{1,2\}$, the morphism $F_{i}$ is given by an
  object $x_{i}$ which is a $G$-torsor $\pi_{i}:P_{i}\to X$ and a
  $G$-equivariant morphism $P_{i}\to V$. Let $\alpha,\alpha':P_{1}\to
  P_{2}$ be morphisms between the objects $x_{1}$ and $x_{2}$ such
  that $\alpha\vert_{\pi_{1}^{-1}(U)}=\alpha'\vert_{\pi_{2}^{-1}(U)}$.
  As $G$ is separated, we have that $\pi_{i}$ is separated. We deduce
  that $\alpha=\alpha'$.
  
  Now we prove the uniqueness of the proposition in the case where
  $\YY=[V/G]$.  Let $X$ be an \'etale atlas of $\XX$. By the previous
  point, we deduce that $\alpha\vert_{X}=\alpha'\vert_{X}$. As
  $\Mor(F_{1},F_{2})$ is a sheaf on $\XX$, we conclude that
  $\alpha=\alpha'$.
  
  For the general case, we reduce to the previous by covering $\YY$ by
  global quotients and then we use that $\Mor(F_{1},F_{2})$ is a sheaf
  on $\XX$.
  
  \textbf{Existence:} It is enough to do it for an \'etale affine
  chart  of $\XX$. By hypothesis, this chart is a disjoint
  union of affine irreducible normal varieties. Hence, we can assume that
  $\XX$ is an affine irreducible normal variety, denoted by $X$.
  Denote by $U$ the scheme $\UU$ open dense in $X$. The morphism
  $F_{1}\circ \iota:U\to \YY$, the $2$-arrow $\beta$ and the
  universal property of the strict fiber product give a morphism
  $f:U\to U'$.  The existence of $\alpha$ is equivalent to the
  existence of a morphism $h:\xymatrix{X \ar@{-->}[r]& X'}$ such that
  $\pi_{1}\circ h =\Id$ and $h\circ \iota=g\circ f$. Denote by
  $\Delta:\YY\to \YY \times \YY$ the diagonal. We can sum up the
  informations in the following diagram :
  \begin{displaymath}
    \xymatrix{
U\ar@/_1pc/[rdd]^{\Id}\ar[rr]^{\iota}\ar[dr]^{f}&
&X\ar@/_1pc/[rdd]^(.2){\Id}
\ar@{-->}[rd]^{\exists\ 
h}&&& \\
&U'\ar[rr]^{g}\ar@{}[rrd]|{\square}\ar[d]&&X'\ar[rr]^{\pi_{2}}\ar[d]^{\pi_{1}}\ar@{}[rrd]|{\square}
&&\YY\ar[d]^{\Delta}\\
&U\ar[rr]^{\iota}&&X\ar[rr]^-{(F_{1}\times F_{2})}&&\YY\times\YY}
  \end{displaymath}
 By definition of the separatedness
  of $\YY$, we have the $\Delta$ is proper. By Lemma 4.2 of
  \cite{LMBca}, we have that $\Delta$ is finite and $X'$ is a scheme.
  We deduce that $\pi_{1}:X' \to X$ is finite. The morphism $g\circ f
  :U\to X'$ is a section of $\pi_{1}$. By Lemma \ref{lem:useful}, we
  deduce a morphism $h:X\to X'$ such that $h\circ \iota=g\circ
  f$. This completes the proof.
  

\end{proof}

\begin{lem}\label{lem:useful} 
  Let $X'$ be a scheme and $X$ be an irreducible normal variety. Let
  $\pi:X'\to X$ a finite morphism. Let $U\hookrightarrow X$ be an open
  dense immersion. Let $s:U\to X'$ be a section of $\pi$. Then the
  section $s$ extends to a section $\widetilde{s}:X\to X'$.
  \end{lem}

\begin{proof}
  Denote by $U_{0}$ the closure of the $s(U)$ in the fiber product
  $U':=U\times_{X}X'$. Denote by $p:U'\to U$ and $q:U'\to X'$ the
  morphisms induced by the fiber product $U'$. Looking at the
  fractional fields, we deduce that the morphisms $s:U\to U_{0}$ and
  $p\vert_{U_{0}}: U_{0}\to U$ are birational morphisms. Denote  by $X_{0}$
  the closure of $U_{0}$ in $X'$. As the morphism $q$ is an open
  embedding, we have that $q\vert_{U_{0}}$ is dominant. We deduce that
  $\pi\vert_{X_{0}}:X_{0}\to X$ is birational and quasi-finite.
  As $X$ be an irreducible normal variety, the Zariski main theorem
  implies that $\pi\vert_{X_{0}}$ is an isomorphism. Its inverse is the
  wanted section of $\pi$. 
\end{proof}

 \section{Action of a Picard stack}\label{sec:action-picard-stack}
 
 In this appendix, we recall the definition of a Picard stack. Then we
 define the action of a Picard stack on a stack which extends the
 definition of Romagny in \cite{Rgasa-2005}. 
In  \cite[Definition 6.1]{Breen-Bitors-1990} Breen defines the notion of a $\mathcal{G}$-torsor over a stack where $\mathcal{G}$ is a Picard stack. Our definition of the action is actually already included in that definition. 
 
 To define the notion of Picard stacks, we do not need the stacks to
 be algebraic.
\begin{defi}[Picard Stacks \cite{SGA4} Exp. XVIII]\label{defi:picard,stacks}
  Let $S$ be a base scheme. A Picard $S$-stack $\mathcal{G}$ is an
  $S$-stack with the following data :
  \begin{itemize}
  \item(Multiplication) a morphism of $S$-stacks:
    \begin{displaymath}
  \xymatrix@R=0pt@1{
    \mathcal{G}\times_S\mathcal{G} \ar[r]^-m & \mathcal{G} \\
    (g_1,g_2) \ar@{|->}[r] & g_1\cdot g_2 \\
    }
\end{displaymath}

\item($2$-Associativity) a 2-arrow $\theta$ implementing the
associativity law:
\begin{equation}\label{eq:associativy.of.picard.stacks}
 \theta_{g_1,g_2,g_3}: (g_1\cdot g_2)\cdot g_3 \Rightarrow g_1\cdot(g_2\cdot g_3)
\end{equation}
\item($2$-Commutativity) a 2-arrow $\tau$ implementing commutativity:
\begin{equation}\label{eq:commutativity.of.picard.stacks}
  \tau_{g_1,g_2}:  g_1\cdot g_2 \Rightarrow g_2\cdot g_1
\end{equation}
\end{itemize}
These data must satisfy the following conditions:
\begin{enumerate}
\item for every chart $U$ and every object $g\in\mathcal{G}(U)$ the map $m_g:\mathcal{G}\to\mathcal{G}$ which multiplies every object by $g$ and every arrow by $\Id_g$ is an isomorphism of stacks.
\item (\textit{Pentagon relation}) For every chart $U$ and 4-tuples of
  objects $g_i\in\mathcal{G}(U)$, we have 
  \begin{equation}
    (\Id_{g_1}\cdot\theta_{g_2,g_3,g_4})\circ\theta_{g_1,g_2\cdot g_3,g_4}\circ(\theta_{g_1,g_2,g_3}\cdot\Id_{g_4})=\theta_{g_1,g_2,g_3\cdot g_4}\circ\theta_{g_1\cdot g_2,g_3,g_4}
  \end{equation}
\item  For every chart $U$ and every object $g\in\mathcal{G}(U)$, we have $\tau_{g,g}=\Id_{g\cdot g}$. 
 \item For every chart $U$ and every objects
 $g_1,g_2\in\mathcal{G}(U)$, we have $\tau_{g_1,g_2}\circ\tau_{g_2,g_1}=\Id_{g_2\cdot g_1}$.
\item (\textit{Hexagon relation}) For every chart $U$ and every triple
  of objects $g_{1},g_{2},g_{3}$ in $\mathcal{G}(U)$, we have
\begin{equation}\label{eq:hexagon}
 \theta_{g_1,g_2,g_3}\circ\tau_{g_3,g_1\cdot g_2}\circ \theta_{g_3,g_1,g_2}=(\Id_{g_1}\cdot \tau_{g_3,g_2})\circ\theta_{g_1,g_3,g_2}\circ(\tau_{g_3,g_1}\cdot\Id_{g_2})
\end{equation}
\end{enumerate}
\end{defi}

\begin{rem}
The Pentagon relation establishes the compatibility law between 2-arrows $\theta$ when
expressing the associativity with 4 objects

The third condition means that every object strictly commutes with itself. 

The last condition states compatibility between the 2-arrow of associativity
and the 2-arrow of commutativity.
\end{rem}

\begin{rem}\label{rem:neutral.element}
  It can be proved, see  \cite[Exp.XVIII,1.4.4]{SGA4},
  that the previous definition is enough to guarantee the existence of
  a neutral element in the group stack. More precisely it is a couple $(e,\epsilon)$ where 
  $e:S\rightarrow\mathcal{G}$ is a section and $\epsilon: e\cdot e
  \Rightarrow e$. A neutral element is unique up to a unique isomorphism. 
\end{rem}

\begin{defi}[Morphisms of Picard stacks \cite{SGA4} Exp. XVIII] Let $(\mathcal{G},\theta,\tau)$ and $(\mathcal{H},\psi,\rho)$ be two Picard $S$-stacks.
  A morphism of Picard $S$-stacks is a morphism of $S$-stacks
  $F:\mathcal{G}\to\mathcal{H}$ with a 2-arrow $\phi_{g_1,g_2} :F(g_1\cdot
  g_2)\Rightarrow F(g_1)\cdot F(g_2)$ for any $g_{1}, g_{2}$ objects
  of $\mathcal{G}$ satisfying the following compatibility conditions:

\begin{itemize}
\item for every chart $U$ and every couple of objects $g_1,g_2\in\mathcal{G}(U)$ we have:
  \begin{equation}
    \label{eq:comp.mor.1}
    \rho_{F(g_1),F(g_2)}\circ \phi_{g_1,g_2}=\phi_{g_2,g_1}\circ F(\tau_{g_1,g_2})
  \end{equation}
\item for every chart $U$ and every triple of objects $g_1,g_2,g_3\in\mathcal{G}(U)$ we have:
  \begin{equation}
    \label{eq:comp.mor.2}
    \phi_{g_1,g_2\cdot g_3}\circ (\Id_{F(g_1)}\cdot \phi_{g_2,g_3})\circ F(\theta_{g_1,g_2,g_3})=\psi_{F(g_1),F(g_2),F(g_3)}\circ (\phi_{g_1,g_2}\cdot \Id_{F(g_3)})\circ \phi_{g_1\cdot g_2,g_3}
  \end{equation}
\end{itemize}

\end{defi}
\begin{rem}\label{rem:neutral.lambda}
  \begin{enumerate}
  \item It should be observed that the morphism $F$ maps the pentagon relation (\resp the hexagon relation) for the Picard stack $\mathcal{G}$ to the pentagon relation (\resp the hexagon relation) for $\mathcal{H}$.
  \item Denote by $(e_{\mathcal{G}},\epsilon_{\mathcal{G}})$ a neutral element of $\mathcal{G}$ and $(e_{\mathcal{H}},\epsilon_{\mathcal{H}})$ a neutral element of $\mathcal{H}$. The couple $(F(e_{\mathcal{G}}),F(\epsilon_{\mathcal{G}})\circ \phi_{e_{\mathcal{G}},e_{\mathcal{G}}}^{-1})$ is a neutral element of $\mathcal{H}$. By Remark \ref{rem:neutral.element} there exists a unique $2$-arrow $\lambda : F(e_{\mathcal{G}})\Rightarrow e_{\mathcal{H}}$ such that $\lambda\circ F(\epsilon_{\mathcal{G}})\circ \phi_{e_{\mathcal{G}},e_{\mathcal{G}}}^{-1})=\epsilon_{\mathcal{H}}\circ\lambda^2$.  
\item It can be useful to notice that given $\alpha: g_1\Rightarrow g_2$ and $\beta: g_3\Rightarrow g_4$ morphisms in $\mathcal{G}(U)$ the following identities involving morphisms holds:
  \begin{displaymath}
    F(\alpha\cdot\beta)=\phi_{g_2,g_4}^{-1}\circ ( F(\alpha)\cdot F(\beta))\circ\phi_{g_1,g_3}
  \end{displaymath}
  \end{enumerate}
\end{rem}
\begin{defi}[Action of a Picard stack]\label{defi:action,picard,stack}
  Let $(\mathcal{G},\tau,\theta)$ be a Picard $S$-stack. Denote by $e$
  the neutral section and by $\epsilon$ the corresponding $2$-arrow.
  Let $\XX$ be a $S$-stack. An action of $\mathcal{G}$ on
  $\mathcal{X}$ is the following data :
  \begin{itemize}
  \item  a morphism of $S$-stack :
\begin{displaymath}
\xymatrix@R=0pt@1{
 \mathcal{G}\times_S\mathcal{X} \ar[rr]^a & & \mathcal{X} \\
    g,x \ar@{|->}[rr] & & g\times x \\
}
\end{displaymath}
\item  a 2-arrow $\eta$:
\begin{equation*}
  \label{eq:neutral.action}
 \eta_x :  e\times x \Rightarrow x
\end{equation*}
\item  a 2-arrow $\sigma$:
\begin{equation*}
  \label{eq:associativity.action}
  \sigma_{g_{1},g_{2},x} :  (g_{1}\cdot g_{2})\times x\Rightarrow g_{1}\times(g_{2}\times x).
\end{equation*}
\end{itemize}
These data must satisfy the following conditions :
\begin{enumerate}
\item(\textit{Pentagon}) For every chart $U$, every  objects
  $g_1,g_{2},g_{3}\in\mathcal{G}(U)$ and every object
  $x\in\mathcal{X}(U)$, we have :  
\begin{equation*}
  \label{eq:pentagon.bis}
  (\Id_{g_1}\times\sigma_{g_2,g_3,x})\circ\sigma_{g_1,g_2\cdot g_3,x}\circ(\theta_{g_1,g_2,g_3}\times\Id_x)=\sigma_{g_1,g_2,g_3\times x}\circ\sigma_{g_1\cdot g_2,g_3,x}
\end{equation*}
\item For any chart $U$ and any object $x\in \XX(U)$, we have 
\begin{equation*}
  \label{eq:compatibility.neutral.action}
  (\Id_e\times\eta_x)\circ\sigma_{e,e,x}=(\epsilon\times\Id_x).
\end{equation*}
\end{enumerate}
\end{defi}


\begin{rem}
  \begin{enumerate}
  \item If the Picard stack is a group-scheme then our definition of
    the action is compatible with the one given by Romagny in
    \cite{Rgasa-2005}.
  \item Let $(\mathcal{G},m,\theta,\tau)$ be a Picard $S$-stack. The
    multiplication $m$ defines an action of $\mathcal{G}$ on itself.
\end{enumerate}
\end{rem}
\begin{prop}
  Let $\mathcal{G}_{1}$ and $\mathcal{G}_{2}$ be two Picard
    $S$-stacks. Let $F:\mathcal{G}_{1}\to \mathcal{G}_{2}$ be a
    morphism of Picard stacks with the $2$-arrow $\phi_{g_1,g_2}: F(g_1\cdot g_2)\Rightarrow F(g_1)\cdot F(g_2)$. Let $\XX$ be a $S$-stack with an
    action of $\mathcal{G}_{2}$ given by $(a,\eta,\sigma)$. Then the morphism $F$ induces a
    natural action of $\mathcal{G}_{1}$ on $\XX$.
\end{prop}
\begin{proof}
 The natural action is given by $(\widetilde{a},\widetilde{\eta},\widetilde{\sigma})$ where we put:
 \begin{itemize}
 \item 
  $\widetilde{a}:=a\circ F$;
\item
for every object $x$ in $\XX$, $\widetilde{\eta}_x:=(\eta_x\circ (\lambda\times\Id_x))$ where $\lambda$ is the $2$-arrow defined in Remark \ref{rem:neutral.lambda};
\item 
for every couple $(g_1,g_2)$ of objects of $\mathcal{G}_1$ and every $x$ object of $\XX$, $\widetilde{\sigma}_{g_1,g_2,x}:=\sigma_{F(g_1),F(g_2),x}\circ (\phi_{g_1,g_2}\times\Id_x)$ .
\end{itemize}  
It is straightforward but tedious to check that the triple so defined satisfies all the properties in Definition \ref{defi:action,picard,stack}.
\end{proof}

We finish this section with a  Proposition about actions on algebraic stacks.  
We refer to Definition 12.1 of \cite{LMBca} for the notion of \'etale
site of a \DM stack. 
\begin{prop}\label{prop:action,morphism,BG}
  Let $\XX$ be a smooth  \DM stack and $G$ a finite abelian group. An action of
  $\mathcal{B}G$ on $\mathcal{X}$ induces a morphism of sheaves of
  groups $j:G\times \mathcal{X} \to \IgenX$ on the \'etale site
  of $\XX$. Moreover, as morphism of stacks, $j$ is \'etale.
 \end{prop}

\begin{proof} We may assume $\XX$ to be irreducible and $d$-dimensional. First we produce a stack morphism 
  $j:\mathcal{X}\times G \rightarrow \IgenX$ and we prove that $j$ is
  \'etale. Denote by $e:\Spec \cc \to \mathcal{B}G$ the neutral
  section.  Denote by $\Delta: \XX \to \XX \times \XX$ the diagonal
  morphism. Denote by $a:\mathcal{B}G\times\XX \to \XX$ the action.
  Using the universal property of the fibered product, we have the
  following $2$-commutative diagram:
\begin{equation}\label{eq:35}
  \xymatrix@1{
    \mathcal{X}\times G \ar@{-->}[dr]^j \ar@/_/[ddr]^{p}  \ar@/^/[rrd]^{p} & & & \\
           &{}\save[]+<0.9cm,-0.6cm>*{\square}\restore I(\mathcal{X})\ar[r]^{\pi_{1}}\ar[d]_{\pi_{2}} & \mathcal{X}\ar[d]^{\Delta}\ar@/^/[ddr]^{\Id\times e} & \\
           & \mathcal{X} \ar[r]^\Delta \ar@/_/[rrd]_{\Id\times e} & \mathcal{X}\times\mathcal{X} & \\
           &              &   &     \mathcal{X}\times \mathcal{B}G \ar[ul]^{\Id\times a} \\
  }
\end{equation} where $p:\XX\times G \to \XX$ is the projection.
The stack morphism $j$ must be unramified since it is a factor of the
\'etale morphism $p:\XX\times G \to \XX$. Since every component of
$I(\XX)$ has dimension at most $d$, the stack morphism $j$ is actually
\'etale and its image is contained in $\IgenX$.

Now, it remains to prove that $j: \mathcal{X} \to \IgenX$ is a
morphism of sheaves of groups on the \'etale site of $\XX$.  The two
upper triangles of Diagram (\ref{eq:35}) are strictly commutative
since $I(\mathcal{X})$ is the strict fibered product. This implies
that $j$ is a morphism of sheaves of sets over $\XX$. Notice that on
the \'etale site, the sheaf $I(\XX)$ is $\IgenX$. 

To finish the proof, we need to show that $j$ is a morphism of sheaves
of groups.
Let us check  the compatibility between the
composition law in $I(\mathcal{X})$ and the multiplication of
$G$. This compatibility follows from the existence of a dashed arrow
such that the upper square in the following diagram is strictly commutative.
\begin{displaymath}
  \xymatrix@1{
\mathcal{X}\times G \times G \ar[ddd] \ar@{-->}[dr] \ar[rrr]^{\Id\times m} & & & \mathcal{X}\times G \ar[dl]^j \ar@/^1pc/[ddl]^{p_2}\\
 &{}\save[]+<1.3cm,-0.6cm>*{\square}\restore I(\mathcal{X})\times_{\mathcal{X}}I(\mathcal{X}) \ar[r]^(.65)c \ar[d] & I(\mathcal{X})\ar[d]^{\pi_2} & \\
 & I(\mathcal{X}) \ar[r]^{\pi_2} & \mathcal{X} &  \\
  \mathcal{X}\times G \ar[ur]_j \ar@/_1pc/[urr]_{p_2} & & & \\
    }
\end{displaymath}
where the stack morphism $c$ is the composition law of the inertia
 stack.
The external square of the diagram above is
 $2$-cartesian and  the stack morphism $\Id \times m :\XX\times G\times G \to
 \XX\times G$ is the identity on $\XX$ and the multiplication in $G$.
 By the universal property of the strict fiber product, we deduce the
 dashed arrow such that the upper square is strictly commutative. This
 ends the proof.

\end{proof}

\section{Stacky version of the Zariski's Main Theorem}
\label{sec:app:Zariski}
Here, we prove a stacky version of  the Zariski's
Main Theorem. We did not
find any reference in the literature for this version.

\begin{thm}[Zariski's Main Theorem for stacks]\label{thm:stack-vers-zariski}
  Let $\mathcal{X}$, $\mathcal{Y}$ be smooth \DM stacks. Let
  $f:\mathcal{X}\to \mathcal{Y}$ be a representable, birational,
  quasi-finite and surjective morphism. Then $f$ is an isomorphism. 
\end{thm}

\begin{proof}
  Let $Y\to\mathcal{Y}$ be an \'etale atlas. Consider the following
  fiber product
  \begin{displaymath}
    \xymatrix@1{
X\ar[r]^-{\overline{f}}\ar[d]\ar@{}[dr]|{\square}&\ar[d] Y \\ 
\mathcal{X}\ar[r]^-{f}& \mathcal{Y}}
  \end{displaymath}
  The morphism $\overline{f}:X\to Y$ is proper, birational, surjective
  and quasi-finite between smooth varieties.  Hence, the Zariski Main
  Theorem (see for example \cite[p.209]{Mumford-red-book}) implies that $\overline{f}$ is an isomorphism. This implies
  that $f:\mathcal{X}\to \mathcal{Y}$ is an isomorphism.
\end{proof}

\bibliographystyle{alpha}
\bibliography{bibliostack}

\newcommand{\etalchar}[1]{$^{#1}$}
\begin{thebibliography}{HKK{\etalchar{+}}03}

\bibitem[$\aleph$CV03]{ACVtbac-2003}
Dan Abramovich, Alessio Corti, and Angelo Vistoli.
\newblock Twisted bundles and admissible covers.
\newblock {\em Comm. Algebra}, 31(8):3547--3618, 2003.
\newblock Special issue in honor of Steven L. Kleiman.

\bibitem[SGA4]{SGA4}
Michael Artin, Alexander Grothendieck, and Jean-Louis Verdier.
\newblock {\em Th{\'e}orie des Topos et Cohomologie {\'E}tale des Sch{\'e}mas}.
\newblock Lecture notes in Mathematics. Springer-Verlag, 1971.

\bibitem[$\aleph$GV08]{AGVgwdms}
Dan Abramovich, Tom Graber, and Angelo Vistoli.
\newblock Gromov-{W}itten theory of {D}eligne-{M}umford stacks.
\newblock {\em Amer. J. Math.}, 130(5):1337--1398, 2008.

\bibitem[AK70]{Altman-Kleiman-Intro-Grothendieck-duality}
Allen Altman and Steven Kleiman.
\newblock {\em Introduction to {G}rothendieck duality theory}.
\newblock Lecture Notes in Mathematics, Vol. 146. Springer-Verlag, Berlin,
  1970.

\bibitem[$\aleph$OV08]{abramovich-2007}
Dan Abramovich, Martin Olsson, and Angelo Vistoli.
\newblock Tame stacks in positive characteristic.
\newblock {\em Ann. Inst. Fourier (Grenoble)}, 58(4):1057--1091, 2008.

\bibitem[Art71]{Ajhr71}
Michael Artin.
\newblock On the joins of hensel rings.
\newblock {\em Advances in Math.}, 7:282--296 (1971), 1971.

\bibitem[$\aleph$V02]{Abramovich-Vistoli-Compactification-stable-maps}
Dan Abramovich and Angelo Vistoli.
\newblock Compactifying the space of stable maps.
\newblock {\em J. Amer. Math. Soc.}, 15(1):27--75 (electronic), 2002.

\bibitem[BB93]{Bal1993}
Andrzej Bia{\l}ynicki-Birula.
\newblock On actions of {${\bf C}\sp *$} on algebraic spaces.
\newblock {\em Ann. Inst. Fourier (Grenoble)}, 43(2):359--364, 1993.

\bibitem[BC07]{BCqh-2007}
Arend Bayer and Charles Cadman.
\newblock {Quantum cohomology of $[\mathbb{C}^N/\mu_r]$}.
\newblock {\em arXiv.math.AG/0705.2160}, 2007.

\bibitem[BCS05]{BCSocdms05}
Lev~A. Borisov, Linda Chen, and Gregory~G. Smith.
\newblock {The orbifold Chow ring of toric Deligne-Mumford stacks}.
\newblock {\em J. Amer. Math. Soc.}, (1):193--215 (electronic), 2005.

\bibitem[Bre90]{Breen-Bitors-1990}
Lawrence Breen.
\newblock Bitorseurs et cohomologie non ab\'elienne.
\newblock In {\em The {G}rothendieck {F}estschrift, {V}ol.\ {I}}, volume~86 of
  {\em Progr. Math.}, pages 401--476. Birkh\"auser Boston, Boston, MA, 1990.

\bibitem[Cad07]{Cstc-2007}
Charles Cadman.
\newblock {Using stacks to impose tangency conditions on curves}.
\newblock {\em Amer. J. Math.}, 129(2):405--427, 2007.

\bibitem[Cox95a]{Chcr95}
David~A. Cox.
\newblock The homogeneous coordinate ring of a toric variety.
\newblock {\em J. Algebraic Geom.}, 4(1):17--50, 1995.

\bibitem[Cox95b]{Cox-functorofsmooth-1995}
David~A. Cox.
\newblock {The functor of a smooth toric variety}.
\newblock {\em Tohoku Math. J. (2)}, 47(2):251--262, 1995.

\bibitem[Cox05]{Cltv2000}
David Cox.
\newblock {Lectures on Toric Varieties}.
\newblock {\em School on Commutative Algebra given in Hanoi,
  www.amherst.edu/$\sim$dacox/}, 2005.

\bibitem[EHKV01]{EKVbgqs-2001}
Dan Edidin, Brendan Hassett, Andrew Kresch, and Angelo Vistoli.
\newblock Brauer groups and quotient stacks.
\newblock {\em Amer. J. Math.}, 123(4):761--777, 2001.

\bibitem[Ful93]{Fitv}
William Fulton.
\newblock {\em Introduction to toric varieties}, volume 131 of {\em Annals of
  Mathematics Studies}.
\newblock Princeton University Press, Princeton, NJ, 1993.
\newblock , The William H. Roever Lectures in Geometry.

\bibitem[Gir71]{Gcna-1971}
Jean Giraud.
\newblock {\em Cohomologie non ab\'elienne}.
\newblock Springer-Verlag, Berlin, 1971.
\newblock Die Grundlehren der mathematischen Wissenschaften, Band 179.

\bibitem[SGA1]{SGA1}
Alexander Grothendieck.
\newblock {\em Rev\^etements \'etales et groupe fondamental. {F}asc. {I}:
  {E}xpos\'es 1 \`a 5}, volume 1960/61 of {\em S\'eminaire de G\'eom\'etrie
  Alg\'ebrique}.
\newblock Institut des Hautes \'Etudes Scientifiques, Paris, 1963.

\bibitem[Gro68]{Ggb68}
Alexander Grothendieck.
\newblock {Le groupe de Brauer. III. Exemples et compl\'ements}.
\newblock In {\em Dix Expos\'es sur la Cohomologie des Sch\'emas}, pages
  88--188. North-Holland, Amsterdam, 1968.

\bibitem[HKK{\etalchar{+}}03]{Hms2003}
Kentaro Hori, Sheldon Katz, Albrecht Klemm, Rahul Pandharipande, Richard
  Thomas, Cumrun Vafa, Ravi Vakil, and Eric Zaslow.
\newblock {\em Mirror symmetry}, volume~1 of {\em Clay Mathematics Monographs}.
\newblock American Mathematical Society, Providence, RI, 2003.
\newblock With a preface by Vafa.

\bibitem[Iwa06]{Icts06}
Isamu Iwanari.
\newblock The category of toric stacks.
\newblock {\em arXiv.math.AG/0610548}, 2006.

\bibitem[Iwa07a]{Iwanari-GeneralizationofCox-2007}
Isamu Iwanari.
\newblock {Generalization of Cox functor}.
\newblock {\em arXiv.org:0705.3524}, 2007.

\bibitem[Iwa07b]{Iwanari-Logarithmic-2007}
Isamu Iwanari.
\newblock Logarithmic geometry, minimal free resolutions and toric algebraic
  stacks.
\newblock {\em arXiv.org:0707.2568}, 2007.

\bibitem[JT07]{Jiang-Tseng-IntegralChowRing-2007}
Yunfeng Jiang and Hsian-Hua Tseng.
\newblock The integral (orbifold) chow ring of toric deligne-mumford stacks.
\newblock {\em ArXiv.org:0707.2972}, 2007.

\bibitem[LMB00]{LMBca}
G{\'e}rard Laumon and Laurent Moret-Bailly.
\newblock {\em Champs alg\'ebriques}, volume~39 of {\em Ergebnisse der
  Mathematik und ihrer Grenzgebiete. 3. Folge. A Series of Modern Surveys in
  Mathematics [Results in Mathematics and Related Areas. 3rd Series. A Series
  of Modern Surveys in Mathematics]}.
\newblock Springer-Verlag, Berlin, 2000.

\bibitem[Mat89]{Mca}
Hideyuki Matsumura.
\newblock {\em Commutative ring theory}, volume~8 of {\em Cambridge Studies in
  Advanced Mathematics}.
\newblock Cambridge University Press, Cambridge, second edition, 1989.
\newblock Translated from the Japanese by M. Reid.

\bibitem[Mil80]{Met}
James~S. Milne.
\newblock {\em {\'E}tale cohomology}, volume~33 of {\em Princeton Mathematical
  Series}.
\newblock Princeton University Press, Princeton, N.J., 1980.

\bibitem[Mum99]{Mumford-red-book}
David Mumford.
\newblock {\em The red book of varieties and schemes}, volume 1358 of {\em
  Lecture Notes in Mathematics}.
\newblock Springer-Verlag, Berlin, expanded edition, 1999.
\newblock Includes the Michigan lectures (1974) on curves and their Jacobians,
  With contributions by Enrico Arbarello.

\bibitem[Ols07]{Osas05}
Martin Olsson.
\newblock {Sheaves on Artin stacks}.
\newblock {\em J. Reine Angew. Math.}, 603:55--112, 2007.

\bibitem[Per08]{Perroni-notetoricDeligne-Mumford-2007}
Fabio Perroni.
\newblock A note on toric {D}eligne-{M}umford stacks.
\newblock {\em Tohoku Math. J. (2)}, 60(3):441--458, 2008.

\bibitem[Rom05]{Rgasa-2005}
Matthieu Romagny.
\newblock Group actions on stacks and applications.
\newblock {\em Michigan Math. J.}, 53(1):209--236, 2005.

\bibitem[Vis89]{Vitas-1989}
Angelo Vistoli.
\newblock Intersection theory on algebraic stacks and on their moduli spaces.
\newblock {\em Invent. Math.}, 97(3):613--670, 1989.

\end{thebibliography}

\end{document}